\documentstyle[amssymb,amstex,graphicx,graphics,12pt,a4wide]{article}
\newtheorem{theorem}{Theorem}

\newtheorem{proposition}{Proposition}
\newtheorem{corollary}{Corollary}
\newtheorem{remark}{Remark}

\newcommand{\Face}{\operatorname{Face}}
\newcommand{\Faces}{\operatorname{Faces}}
\newcommand{\Facet}{\operatorname{Facet}}
\newcommand{\MaxFaces}{\operatorname{MaxFaces}}
\newcommand{\MaxFacets}{\operatorname{MaxFacets}}
\newcommand{\Cells}{\operatorname{Cells}}
\newcommand{\TypicalCell}{\operatorname{TypicalCell}}
\newcommand{\TKF}{\operatorname{Typical\ k-Face}}
\newcommand{\Vol}{\operatorname{Vol}}

\newcommand{\const}{\operatorname{const}}
\newcommand{\Var}{\operatorname{Var}}
\newcommand{\Ei}{\operatorname{Ei}}
\newcommand{\PLT}{\operatorname{PLT}}
\newcommand{\PHT}{\operatorname{PHT}}
\newcommand{\PVT}{\operatorname{PVT}}
\newcommand{\Cell}{\operatorname{Cell}}
\newcommand{\bd}{\operatorname{bd}}
\newcommand{\diam}{\operatorname{diam}}

\def\Comment#1{
\marginpar{$\bullet$\quad{\tiny #1}}}

\begin{document}

\title{Typical Geometry, Second-Order Properties and Central Limit Theory for Iteration Stable Tessellations}
\author{Tomasz Schreiber\footnote{The first author was supported by the Polish Minister of Science and Higher Education grant N N201 385234 (2008-2010)},\\
Faculty of Mathematics and Computer Science,\\
Nicolaus Copernicus University, Toru\'n, Poland,\\
\textit{e-mail: tomeks at mat.umk.pl};\\
Christoph Th\"ale\footnote{The second author was supported by the Swiss National Science Foundation grant SNF PP002-114715/1.},\\
Department of Mathematics,\\
University of Fribourg, Fribourg, Switzerland,\\
\textit{e-mail: Christoph.Thaele at unifr.ch}}
\date{}
\maketitle

\begin{abstract}
Since the seminal work \cite{NW03,NW05} the iteration stable (STIT) tessellations have attracted considerable interest in stochastic geometry as a natural and flexible yet analytically tractable model for hierarchical spatial cell-splitting and crack-formation processes.  The purpose of this paper is to describe large scale asymptotic geometry of STIT tessellations in ${\Bbb R}^d$ and more generally that of non-stationary iteration infinitely divisible tessellations. We study several aspects of the typical first-order geometry of such tessellations resorting to martingale techniques as providing a direct link between the typical characteristics of STIT tessellations and those of suitable mixtures of Poisson hyperplane tessellations. Further, we also consider second-order properties of STIT and iteration infinitely divisible tessellations, such as the variance of the total surface area of cell boundaries inside a convex observation window. Our techniques, relying on martingale theory and tools from integral geometry, allow us to give explicit and asymptotic formulae. Based on these results, we establish a functional central limit theorem for the length/surface increment processes induced by STIT tessellations. We conclude a central limit theorem for total edge length/facet surface, with normal limit distribution in the planar case and non-normal ones in all higher dimensions.
\end{abstract}
\begin{flushleft}\footnotesize
\end{flushleft}

\section{Introduction}

\subsection{General Introduction}

Infinite divisibility or stochastic stability of a random object under a certain operation is one of the most fundamental concepts in probability theory, prominent examples including the classical theory of infinite divisible and stable distributions with their applications around the central limit theorem, max-stable distributions studied in extreme value theory or union infinitely divisible random sets studied in the classical theory of random closed sets \cite{Mat}.\\
In the present paper we will deal with non-stationary iteration infinitely divisible random tessellations of the $d$-dimensional Euclidean space and more specifically with stationary random tessellations that are \textbf{st}able under the operation of \textbf{it}eration -- called STIT tessellations for short. In the stationary case, the purely mathematical motivation for this type of random tessellations goes back to R. Ambartzumian in the 80thies. The principle of iteration of tessellations can roughly be explained as follows: Take a random primary or frame tessellation and associate with each of its cells an independent copy of the tessellation itself, a component tessellation, which is also independent of the primary tessellation as well. Now make in each cell a local superposition of the primary tessellation and the associated tessellation, whereupon scale the resulting random tessellation by a factor $2$ in order to ensure that the mean surface measure of cell boundaries stays constant. The described operation can now be repeatedly applied and we obtain in this way a sequence of random tessellations. It can be shown that this sequence converges to a random limit tessellation and that this tessellation must be stable under iterations -- a STIT tessellation -- in the translation-invariant set-up. In the general case, the resulting tessellation is iteration infinitely divisible in any finite volume, where this property has the same relation to iteration stable random tessellations as infinitely divisible random variables have to stable ones.\\ Starting with \cite{NW03}, STIT tessellations and their theoretical framework were formally introduced in \cite{NW05} by W. Nagel and V. Weiss. It was the same research group who discovered a first basic technique for studying mean values and even some distributions related to the geometry of STIT tessellations by writing certain balance equations based on the stochastic stability of the tessellation, see \cite{NW06}, \cite{NW08} and \cite{TW} for the mean values as well as \cite{MNW07}, \cite{M09} and \cite{T09} for distributional results. In \cite{MNW} and \cite{MNW2} a new aspect was introduced into the theory, namely a tessellation-valued random Markov process on the positive real half-axis with the property that at each time the law of of the tessellation is stable under iteration. This process sheds light on the hierarchical and temporal structure of STIT tessellations and leads to a random process of cell divisions in any finite volume and moreover in the whole space. This point of view can be exploited to establish further results on STIT tessellations and in fact we will also make use of it in the present paper. In addition, the finite volume Markovian construction provides a link to the class of more general and non-stationary iteration infinitely divisible random tessellations.\\ In next Subsection \ref{secSTIT} we outline important facts about STIT tessellations and introduce the concept of iteration infinitely divisible random tessellations. Afterwards, in Subsection \ref{secOVERVIEW} the plan of the paper as well as a survey of our main results are presented.

\subsection{Iteration Infinitely Divisible and STIT Tessellations}\label{secSTIT}

A tessellation of ${\Bbb R}^d$ is a locally finite partition of the space into compact convex polytopes, the cells of the tessellation. One can regard a tessellation either as a collection of its cells or as the closed set of their boundaries. We will mostly follow the second mentioned path and denote by $\Cells(Y)$ the set of cells of the tessellation $Y$ (by the Jordan--Sch\"onflies theorem the correspondence between $Y$ and $\Cells(Y)$ is one-to-one). Thus, a random tessellation can be regarded as a special random closed set in the sense of \cite{SW}. In particular, this imposes the usual Fell topology and the corresponding Borel measurable structure on the
family of tessellations, see ibidem. A random tessellation $Y$
(regarded as a random closed set in ${\Bbb R}^d$) is called \textit{stationary} if its distribution does not change upon actions
of translations. Analogously a random tessellation is called \textit{isotropic} if its distribution is invariant under the action of $SO(d)$.\\ Whenever we have two random tessellations $Y_1$ and $Y_2$ of ${\Bbb R}^d$ we can define their \textit{iteration/nesting}. To do so, we associate to each cell $c \in \Cells(Y_1)$ an independent version $Y_2(c)$ of $Y_2$ and we assume furthermore the family $\{Y_2(c):c\in \Cells(Y_1)\}$ to be independent of $Y_1$. Then we can define the iteration of $Y_1$ with $Y_2$ by $$ Y_1\boxplus Y_2:= Y_1\boxplus\{Y_2(c):c\in \Cells(Y_1)\}:=Y_1\cup\bigcup_{c\in \Cells(Y_1)}(Y_2(c)\cap c),$$ i.e. we take the local superposition of $Y_2$ and the family $\{Y_2(c):c\in \Cells(Y_1)\}$ inside the cells of $Y_1$. It was shown in \cite{MNW2} that $Y_1\boxplus Y_2$ is a stationary random tessellation as soon $Y_1$ and $Y_2$ are. A stationary random tessellation $Y$ is called \textit{stable under iteration} or \textit{STIT} for short iff
\begin{equation}\label{ITSTAB}
 m \underbrace{(Y \boxplus \ldots \boxplus Y)}_m\overset{D}{=} Y,\; m=2,3,\ldots,
\end{equation}
where $\overset{D}{=}$ stands for equality in distribution, i.e. if its distribution does not change
 under rescaled iteration. In fact, using the uniqueness results, see Theorem 3 and Corollary 2 in
 \cite{NW05}, it is easy to see that it is enough to take one fixed $m > 1$ in (\ref{ITSTAB}).\\ To proceed, let us be given a constant $0<t<\infty$ and a probability measure $\cal R$
on the unit sphere ${\cal S}_{d-1}$ usually identified with the induced
distribution of orthogonal hyperplanes on the space ${\cal H}_0$ of
$(d-1)$-dimensional linear hyperplanes in ${\Bbb R}^d,$ also denoted by
${\cal R}$ in the sequel for notational simplicity. Define the measure $\Lambda$ on the space $\cal H$ of affine hyperplanes
in ${\Bbb R}^d$ as the product measure 
\begin{equation}\label{LADEF}
 \Lambda:=\ell_+ \otimes {\cal R}
\end{equation}
of $\ell_+$ standing for the Lebesgue measure on the positive real half-axis $(0,\infty),$
and of $\cal R$, where a pair $(r,u) \in (0,\infty) \times {\cal S}_{d-1}$ 
is identified with the hyperplane $\{ x \in {\Bbb R}^d,\; \langle x, u \rangle = r \}.$ Throughout this paper we always require that the support of $\cal R$ spans the whole space, i.e. $$\text{span}(\text{supp}({\cal R}))={\Bbb R}^d.$$ Assume now that we are given a stationary random tessellation $Y$ with surface intensity $t$ (i.e. the mean surface area of cell boundaries per unit volume equals $t$) and directional distribution $\cal R$ (i.e. the distribution of the normal direction of the face containing the typical point is given by $\cal R$) and define the sequence $({\cal I}_n(Y))$ by $${\cal I}_1(Y):=2(Y\boxplus Y),\ \ \ \ {\cal I}_{n}(Y):= \frac{n}{n-1}{\cal I}_{n-1}(Y) \boxplus n Y = n \underbrace{(Y \boxplus \ldots \boxplus Y)}_{n},\ n\geq 2.$$
It was shown in \cite[Thm 3]{NW05} that ${\cal I}_n(Y)$ converges in law, as $n \to \infty,$ to a stationary random limit tessellation $Y(t\Lambda)$ {\it uniquely determined} only by $t\Lambda.$ This tessellation is easily shown to be stable under iterations and is called the STIT tessellation with parameters $t$ and $\Lambda$. Without confusion we will write $Y(t)$ instead of $Y(t\Lambda)$, whenever the measure $\Lambda$ or $\cal R$ is fixed.
 \begin{figure}[t]
 \begin{center}
 \includegraphics[width=7cm]{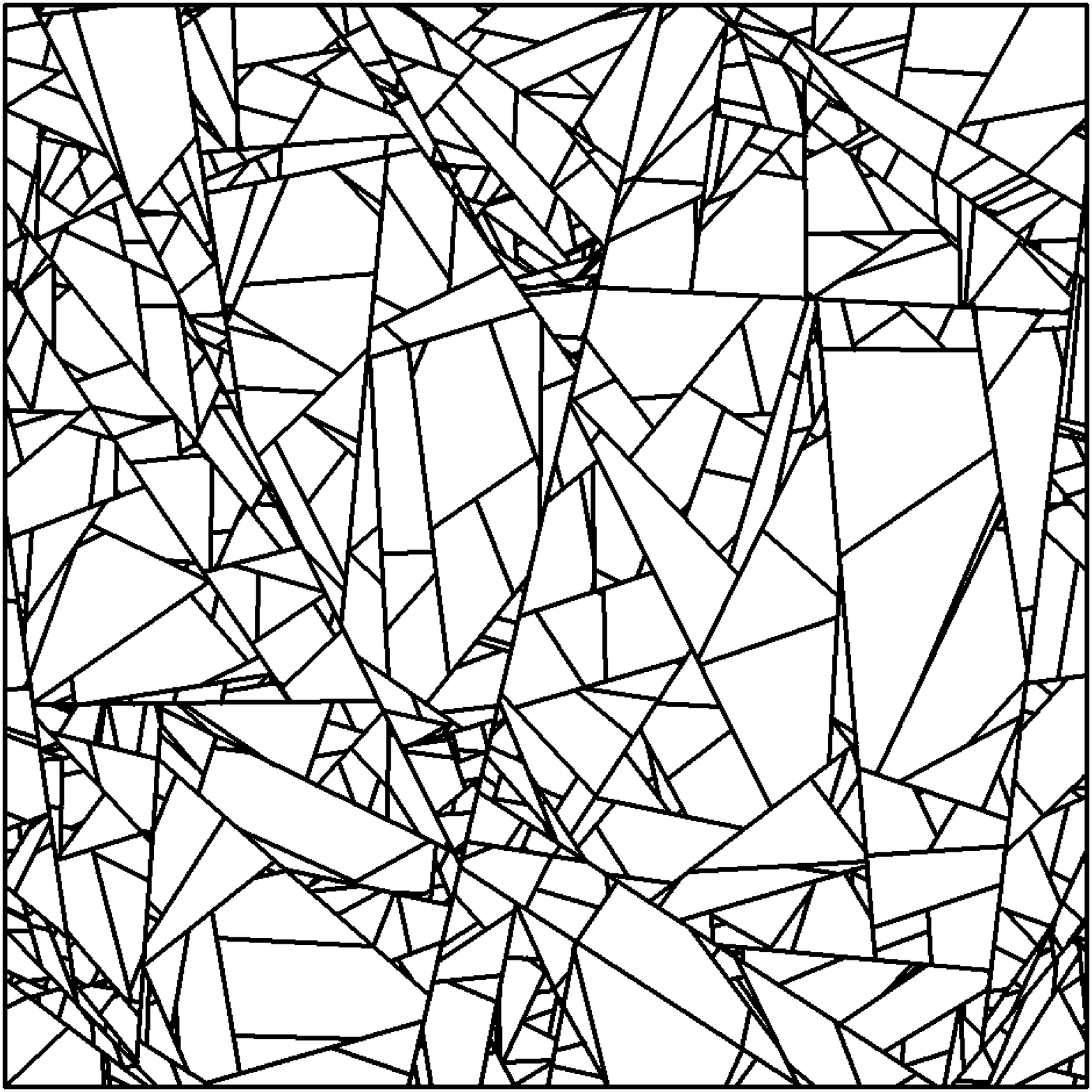}
 \includegraphics[width=7cm]{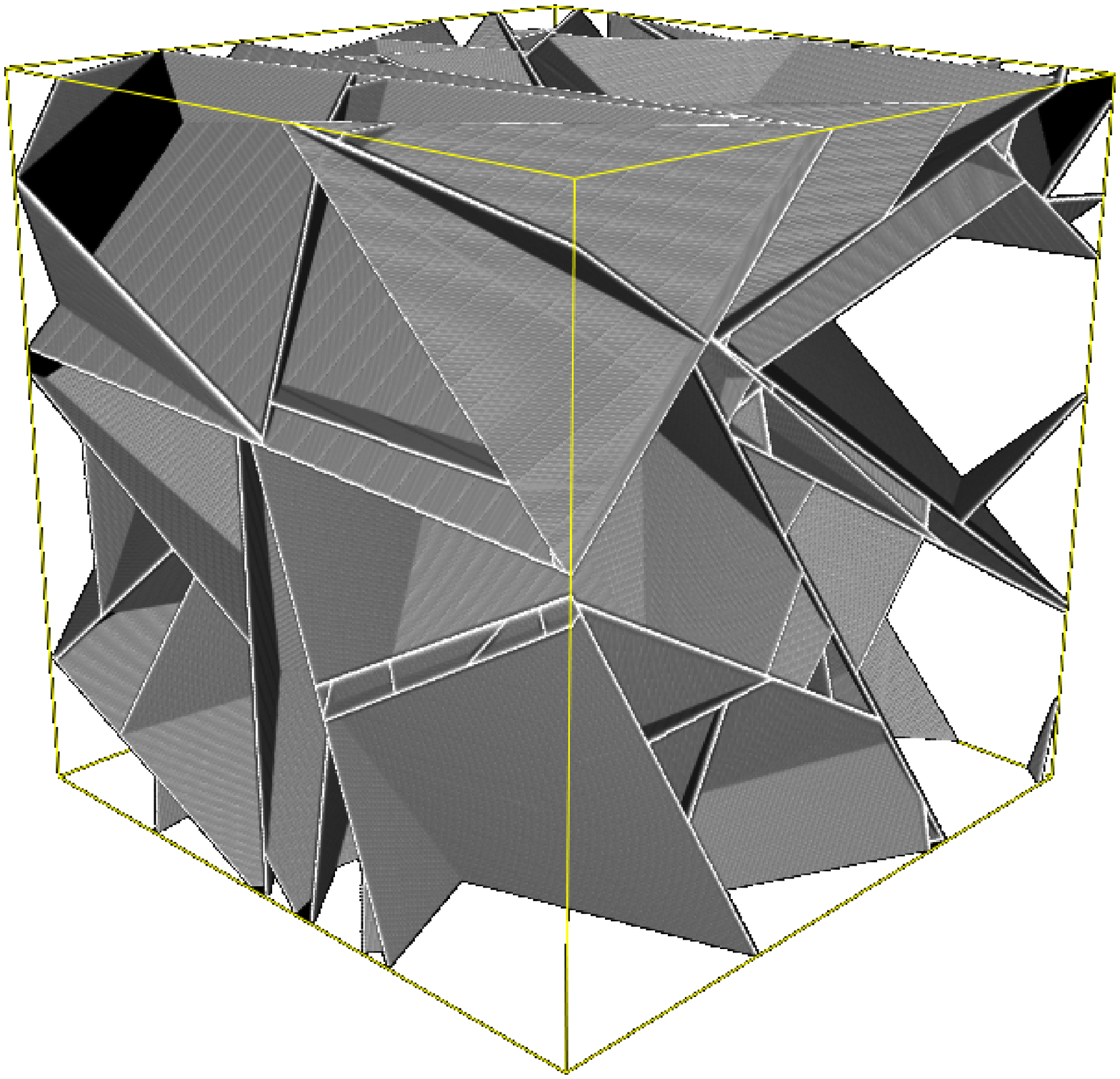}
 \caption{Realizations of a planar and a spatial stationary and isotropic STIT tessellation (kindly provided by Joachim Ohser and Claudia Redenbach)}\label{Fig1}
 \end{center}
\end{figure}
It is a crucial feature
of $Y(t)$ that it admits a very natural and intuitive {\it explicit construction}. For a restriction $Y(t,W)$
of $Y(t)$ to a compact convex window $W\subset{\Bbb R}^d$ this construction can be informally be described as follows
(the reader is referred to \cite{NW05} for full details). Assign to the window $W$ an exponentially distributed
random lifetime with parameter $\Lambda([W])$ where $[W] := \{ H \in {\cal H},\; H \cap W \neq \emptyset \}$
stands for the family of all hyperplanes hitting $W.$ 
 Upon expiry of its lifetime, the cell $W$ dies and splits into two sub-cells $W^+$ and $W^-$ separated by a 
hyperplane in $[W]$ chosen according to the law $\Lambda(\cdot)/\Lambda([W]).$ The resulting 
new cells $W^+$ and $W^-$ are again assigned independent exponential lifetimes with respective
parameters $\Lambda([W^+])$ and $\Lambda([W^-])$ (whence smaller cells live stochastically longer)
and the entire construction continues recursively, until the deterministic time threshold $t$ is reached. 
The cell-separating $(d-1)$-dimensional facets arising in subsequent splits are usually referred to as 
I-faces (or I-segments for $d=2$ as assuming shapes similar to the letter {\it I}).   
The described process of recursive cell divisions is called the Mecke-Nagel-Weiss- or (MNW)-construction in the sequel
and the resulting random tessellation created inside $W$ is denoted by $Y(t,W)$ as mentioned above,
whereas the collection of all I-facets or I-segments is denoted by $\MaxFacets(Y(t,W))$.
Moreover, we write $\MaxFaces_k(Y(t,W))$ for the collection of $k$-dimensional I-faces of $Y(t,W)$,
where by a $k$-dimensional I-face we mean the maximal union of connected and $k$-coplanar $k$-dimensional
faces. In fact, $k$-dimensional I-faces of $Y(t,W)$ can also be alternatively defined as the $k$-faces of I-facets.
 It was shown in \cite{NW05} that the law of $Y(t,W)$ is consistent in $W,$ i.e.
$Y(t,W) \cap V \overset{D}{=} Y(t,V)$ for convex $V \subset W$ and thus $Y(t,W)$ can be extended to random
tessellation $Y(t)=Y(t\Lambda)$ on the whole space, which is then proved (cf. \cite{NW05}) to coincide with the limit tessellation
$Y(t)$ considered above, as the notation suggests. 
Again, the sets of all I-facets and I-faces of $Y(t)$ are denoted by $\MaxFacets(Y(t))$ and $\MaxFaces_k(Y(t))$
($0\leq k\leq d-2$), respectively. 
The stationary random tessellation $Y(t)$ is additionally isotropic if and only if $\cal R$ is the uniform distribution on ${\cal S}_{d-1}$. In this case we will write $\Lambda_{iso}$ or $\Lambda_{iso}^{{\Bbb R}^d}$ for the invariant measure on the space of hyperplanes $\cal H$.\\ A simple yet crucial observation is that even though only translation-invariant measures $\Lambda$ of the form (\ref{LADEF}) show up in the limiting STIT tessellations, the MNW-construction can be
carried out with arbitrary non-atomic and locally finite driving measure $\Lambda$ (i.e. $\Lambda([W]) < \infty$
for $W$ bounded) on ${\cal H}$ also leading to a consistent family $Y(t,W)$ and eventually, by
extension, yielding $Y(t).$ Many of our theorems below will be stated in this general context. It should
be emphasised though that such tessellations are no more iteration stable (STIT). However,
they have the general property of being {\it iteration infinitely divisible}, as they can
be readily checked to arise as $m$-fold iterations of $Y(t/m)$ for each $m \geq 2$ in all
finite volumes. Formally, this means that $$Y(t/m,W)^{\boxplus m} \overset{D}{=} Y(t,W)$$
for all compact convex sets $W \subseteq {\Bbb R}^d,$ which follows directly by the MNW-construction as yielding $$Y(s,W) \boxplus Y(u,W) \overset{D}{=} Y(s+u,W).$$ It is more
than natural to expect that also $Y(t) = Y(t/m)^{\boxplus m}$ for all $m$ in the
whole ${\Bbb R}^d$ which should be easily provable by adopting the theory developed
in \cite{MNW2}, thus even better justifying the term {\it iteration infinitely divisible}
in our context, yet this falls beyond the scope of the present work.\\ The STIT tessellations enjoy a number of remarkable properties. Below, we recall some of them which will be
of importance for our further argumentation. All proofs can be found in \cite{NW03} and \cite{NW05}.
\begin{itemize}
	\item[(a)] STIT tessellations have Poisson typical cells, i.e. under stationary $\Lambda$ the interior of the typical cell of $Y(t)$ - denoted by $\TypicalCell(Y(t))$ - has the same distribution as the interior of the typical cell of a stationary Poisson hyperplane (or line or plane) tessellation with the same surface intensity $t>0$ and the same directional distribution $\cal R, $ see Lemma 3 in \cite{NW05}. The difference between these two types of tessellations arises from the mutual arrangement of the cells, see Fig. \ref{Fig1}.
	\item[(b)] Sectional STIT tessellations are STIT tessellations as well. In particular, 
    let $E_k$ be a $k$-dimensional plane in ${\Bbb R}^d$, $1\leq k\leq d-1$, and let $Y(t)$ be a stationary and isotropic random STIT tessellation with parameter $t>0$. Then $Y(t)\cap E_k$ has the same distribution as the stationary and isotropic STIT tessellation with parameter $\lambda_k t$ constructed inside the plane $E_k$,
    i.e. $$Y(t \Lambda_{iso}^{{\Bbb R}^d},{\Bbb R}^d)\cap E_k\overset{D}{=}Y(\lambda_k t \Lambda_{iso}^{E_k},E_k)$$ with,  cf. \cite[(3.29T)]{MRE}
\begin{equation}
\lambda_k={\Gamma\left({k+1\over 2}\right)\Gamma\left({d\over 2}\right)\over\Gamma\left({k\over 2}\right)\Gamma\left({d+1\over 2}\right)}\label{LAMBDAK}
\end{equation}
by abusing the above introduced notation, where we have assumed the window to be compact. Moreover, if $k=1$ the intersection of $E_1$ with $Y(t)$ induces a motion-invariant Poisson process with intensity $$\lambda_1t={\Gamma\left({d\over 2}\right)\over\sqrt{\pi}\Gamma\left({d+1\over 2}\right)}t$$ on the line $E_1$. In view of the above considerations we see that the motion-invariant STIT tessellations in dimension $1$ are the motion-invariant Poisson processes (note that in dimension $1$ motion-invariance is the same as translation-invariance). Thus, STIT tessellations can be seen as another generalization of the $1$-dimensional Poisson process to higher dimensions beside other models, see \cite{MM}.
	\item[(c)] STIT tessellations have the following scaling property: $$tY(t)\overset{D}{=}Y(1),$$ i.e. the 
 tessellation $Y(t)$ of surface intensity $t$ upon rescaling by factor $t$ has the same distribution as $Y(1)$, the STIT tessellation with surface intensity $1$.
\end{itemize}
We will use in the paper the same notation for iteration infinitely divisible random tessellation and STIT tessellations hopefully without confusion, mentioning in any case if we are working with a general locally finite non-atomic measure $\Lambda$, a translation-invariant one or even with $\Lambda_{iso}$.\\ We will also make use of the following notation in the paper:
\begin{itemize}
	\item $B_R=B_R^d(o)$ is the $d$-dimensional ball around the origin with radius $R>0$.
	\item The $k$-dimensional volume measure will be denoted by $\Vol_k$.
    \item $\kappa_j := \Vol_j(B_1^j)$ is the volume of the $j$-dimensional unit ball, $j\kappa_j$ its surface area.	
    \item The uniform probability measure on the unit sphere ${\cal S}_{d-1}$ in ${\Bbb R}^d$ (normalized spherical surface measure) is denoted by $\nu_{d-1}$.
\end{itemize}

\subsection{Short Overview of the Paper and the Results}\label{secOVERVIEW}

As already remarked above, a STIT random tessellation induces a continuous time tes\-sel\-la\-tion-valued Markov process on the positive real half-axis. For such processes, the notion of a generator is available. Using standard theory of Markov processes, several martingales associated with the tessellations under consideration can be constructed. This fundamental observation will be exploited in Section \ref{secMARTINGALE} and most of our results will be based on it.\\ In Section \ref{secFIRST} we mainly deal with first-order properties and the typical geometry of iteration infinitely divisible or more specifically stationary and stationary and isotropic STIT tessellations. We generalize some properties known for STIT tessellations to the non-stationary case and derive explicit formulas for mean values and higher moments of $k$-dimensional I-faces in the isotropic case (Sections \ref{sectypifaces} and \ref{secHIGERMOMENTS}). In particular we show that the distribution of the (time-marked) typical $k$-dimensional I-face is a mixture of suitably (marked and) rescaled $k$-dimensional Poisson cells. As a main tool we use the martingale techniques introduced before and we compare the tessellations with certain mixtures of Poisson hyperplane tessellations (Section \ref{seclowerdimfaces}). Typical geometry beyond the isotropic regime is explored in Section \ref{SECBEHINDISO}, where for example the conditional distribution of the typical $k$-dimensional I-segment given its birth time for all space dimensions is determined and some global and local mean value formulas for the stationary as well as for the non-stationary case are derived.\\ Second-order parameters of iteration infinitely divisible and STIT tessellations are the contents of Section \ref{secSECOND}. Based on a specialization of our martingale technique (Section \ref{secSECONDMARTINGALETOOLS}), we calculate the variance of a general face-functional and as a special case we find the variance of the total surface area of cell boundaries in a bounded convex window. The resulting integral expression can be explicitly evaluated in the stationary and isotropic case by applying the affine Blaschke-Petkantschin formula (Section \ref{subsecVAR}). For the particular case of space dimension $3$, an exact formula without further integrals is presented as well as general variance asymptotics relying on techniques from integral geometry. Certain chord-power integrals will reflect the influence of the geometry of the observation window to this variance asymptotics. The purpose of Section \ref{secALTERNATIVE} is to compare our martingale approach with a recent elegant second-order theory for stationary and isotropic planar STIT tessellations independently developed by W. Weiss, J. Ohser and W. Nagel. As a by-product, we derive there an explicit expression for the pair-correlation function of the random surface area measure in general space dimensions. For space dimensions $2$ and $3$ we also calculate the so-called $K$-function.\\ Following the afore-mentioned first- and second-order theory we turn in Section \ref{secCLT} to the central limit problem for STIT tessellations. In Section \ref{secBROWNIANCONV} we show the convergence of the surface/length increment process wrt. to some initial time instant associated with STIT tessellation to a time changed Wiener process, the proof of which is based on the ergodicity property of STIT tessellations and the functional central limit theorem for square-integrable martingales. It is interesting to see that in terms of the time incremental MNW-construction of STIT tessellations our central limit theorem is more reliant on independencies arising in construction time than those of spatial nature. Even more interestingly, the 'big-bang' phase near time zero turns out to play a crucial r\^ole. In Section \ref{CLTplane} we establish normal convergence of the total edge length in the planar case and in Section \ref{CLTspace} it is shown that non-normal limit distributions appear in the limit for space dimensions greater or equal $3$. The argument goes by using our variance calculations from Section \ref{secSECOND} and the observation that for space dimensions $\geq 3$ the initial 'big-bang' phase of the MNW-construction brings a non-negligible contribution to the variance of the surface increment process, whereas in dimension $2$ this contribution is negligible but already with an extremely slow rate of decay.\\ The aim of the last Section (Section \ref{secCOMPARISON}) is to put our results in a more general context by comparing them with the results for Poisson hyperplane tessellations and for Poisson-Voronoi tessellations with the same surface intensity.\\ In all cases, our general statements are illustrated by examples in space dimensions $d=2$ and $d=3$, sometimes leading to already known formulas but often extending them at least in some aspects.\\ Our results could -- beside their intrinsic mathematical motivation -- be of potential interest for applications for example in context of questions concerning statistical model fitting of random tessellations to real data. Central limit theorems could be a basis for statistical inference of tessellation models and related functionals. Asymptotic confidence intervals and tests with respect to mean values can be derived from them, since we can make the first- and second-order moments explicitly available.

\section{Martingales Associated with Iteration Infinitely Divisible or STIT Tessellations}\label{secMARTINGALE}
 The finite volume continuous time incremental MNW-construction of iteration infinitely divisible
 random tessellations or more specially stationary STIT tessellations, as discussed in Section \ref{secSTIT} above,
 clearly enjoys the Markov property in time parameter, whence natural martingales arise which
 will be of crucial importance for our further considerations.
 To discuss these processes we need some additional terminology and notational conventions. 
 First, we fix a compact convex window $W$ and a general diffuse (non-atomic) and locally finite measure $\Lambda$
 on ${\cal H}.$ Next, whenever a new facet $f$ lying on a hyperplane $H$ is born
 in the course of the Mecke-Nagel-Weiss- (MNW)-construction, that is to say a cell splits, we declare that
 its {\it initial/birth point} is the point (vertex)
 with lowermost first coordinate, to be denoted $\iota(f)$ in the sequel (in fact any other deterministic
 and measurable algorithm of distinguishing one point of $f$ would do as well in this context).
 Moreover, for a given tessellation $Y$ (usually taken to be $Y(t,W)$ for some $t > 0$) we also introduce the
 notation $\iota(H;Y)$ for $H$ \--- a hyperplane in $[W],$ with the following meaning
 \begin{itemize}
  \item  The hyperplane $H$ is tessellated by the intersection with $Y,$ and thus split
            into a number of $(d-1)$-dimensional polyhedra $f_1,f_2,\ldots$ which could potentially
            become new facets for the tessellation, were its construction continued.
 \item   We put $\iota(H;Y) := \{ \iota(f_1),\iota(f_2),\ldots \}.$
\end{itemize}
Whenever $x \in \iota(H;Y),$ we define $\Facet(x,H|Y)$ to be the new facet to be added to
$Y$ should a facet birth (cell split) occur at $x$ on $H$ in $Y.$ Moreover, we denote by $\Cell(x,H|Y)$
the cell split by the birth of $\Facet(x,H|Y)$ and we write $\Cell^+(x,H|Y)$ and $\Cell^-(x,H|Y)$ for
the two sub-cells into which $\Cell(x,H|Y)$ divides, lying on the positive and negative sides of $H$,
respectively.\\ With this notation, it is easily seen that for a fixed measure $\Lambda$, 
$(Y(t,W))_{t\geq 0}$ is a pure jump Markov process with values in the space of tessellations of $W$ and
with the property that
\begin{equation}\label{INFINEQ1}
  d{\Bbb P}\left(Y(t+dt)=Y \cup \Facet(x;H|Y)|Y(t)=Y\right) = 
   {\bf 1}[x \in \iota(H;Y)]\Lambda(dH)dt, 
\end{equation}
with $H \in [W]$, and
\begin{eqnarray}
   \nonumber {\Bbb P}\left(Y(t+dt) = Y | Y(t) = Y\right) &=& 1 - \left( \int_{[W]} |\iota(H;Y)| \Lambda(dH) \right) dt\\
   &=& 1 - \left( \sum_{c \in \Cells(Y)} \Lambda([c]) \right) dt \label{INFINEQ2}
\end{eqnarray}
with $|\iota(H;Y)|$   standing for the cardinality of $\iota(H;Y).$
Indeed, this is because, conditionally on $Y(t) = Y,$ during the period $(t,t+dt]$ of the MNW-construction
we have, as discussed in Subsection \ref{secSTIT},
\begin{itemize}
 \item 
  For each cell $c \in \Cells(Y)$ the probability that it undergoes a split is $\Lambda([c]) dt,$
  moreover the probability that two or more cells split is $o(dt),$ whence the probability that no split
  occurs is $ 1 - \left( \sum_{c \in \Cells(Y)} \Lambda([c]) \right) dt$ as in (\ref{INFINEQ2}).
  Observing that $|\iota(H;Y)|$ coincides with the number of cells in $Y$ that $H$ intersects,
  we get in addition $\sum_{c \in \Cells(Y)} \Lambda([c]) = \int_{[W]} |\iota(H;Y)| \Lambda(dH)$ which yields 
  the remaining equality in (\ref{INFINEQ2}).
  \item Should a cell $c \in \Cells(Y)$ split, the splitting hyperplane $H$ is chosen according to the law
          ${\bf 1}[{H \in [c]}](\cdot)/\Lambda([c]),$ whence the probability of observing a split of $c$
         induced by $H$ during the time period $(t,t+dt]$ is just ${\bf 1}[{H \in [c]}](dH) dt.$ Now, 
         having $x \in \iota(H;Y)$ is by definition equivalent to there being a cell $c = c(x,H) \in \Cells(Y)$ 
         with $H \in [c]$ and $x = \iota(c \cap H).$ For so chosen  $c$ we see that on the event $Y(t) = Y,$ having
         $Y(t+dt) = Y \cup \Facet(x;H|Y)$ is equivalent to having $c$ split by $H$ during $(t,t+dt].$ As noted
         above, the latter happens with probability $\Lambda(dH) dt$ whence (\ref{INFINEQ1}) follows.
\end{itemize} 
It should be emphasised at this point that the {\it differential} notation in (\ref{INFINEQ1}) and (\ref{INFINEQ2}) 
employing the symbols $d{\Bbb P}$ and $dt$ is widely accepted in the theory of pure jump continuous
time Markov processes and makes perfect formal sense as a commonly recognised abbreviation for
the usual description in terms of waiting times, with $d{\Bbb P}(Y(t+dt) = Y'|Y(t)=Y) = a(Y,Y') dt$
understood as '{\it while in state $Y,$ wait an exponential time with parameter $a(Y,Y')$ and then
jump to $Y'$ unless some other jump has occurred prior to that}'. We refer the reader to Chapter
15 in \cite{BREI} and especially to Section 15.6 there, where this construction is formalised with
full mathematical rigour in terms of waiting times as noted above. As readily verified, specialising
the generic waiting-time construction to the case of (\ref{INFINEQ1})) and (\ref{INFINEQ2}) yields
precisely the standard MNW construction.\\
Using (\ref{INFINEQ1}) and (\ref{INFINEQ2}) we conclude by general theory of Markov processes
and their infinitesimal generators, see Chapter 1 in \cite{LIG} or Chapter 15 and especially Sections
15.4 (Def. 15.21) and 15.6 in \cite{BREI} specialised for the pure jump case, the generator for
$(Y(t,W))_{t\geq 0}$ is ${\Bbb L} := {\Bbb L}_{\Lambda; W}$
\begin{equation}\label{GEN}
 {\Bbb L}F(Y) = \int_{[W]} \sum_{x \in \iota(H;Y)} [F(Y \cup \Facet(x;H|Y)) - F(Y)] \Lambda(dH)
\end{equation}
for all $F$ bounded and measurable on space of tessellations of $W.$
Consequently, again by standard theory as given in Lemma 5.1 Appendix 1 Sec. 5 in \cite{KipLan},
see also Section 1.5 in \cite{LIG}, or alternatively by a direct check  straightforward
in the present set-up, we readily see that for $F$ bounded and measurable the stochastic process
\begin{equation}\label{MART}
 F(Y(t,W)) - \int_0^t {\Bbb L}F(Y(s,W)) ds
\end{equation}
is a martingale with respect to the filtration $\Im_t$ generated by $(Y(s,W))_{0 \leq s \leq t}.$
More generally, for bounded measurable $G = G(Y,t),$ considering the {\it time-augmented} Markov process 
$(Y(t,W),t)_{t \geq 0}$ and applying usual theory, see again Lemma 5.1 in Appendix 1 Sec. 5
in \cite{KipLan}, or simply by performing a direct check, we see that
\begin{equation}\label{MART2}
 G(Y(t,W),t) - \int_0^t \left( [{\Bbb L}G(\cdot,s)](Y(s,W)) + \frac{\partial}{\partial s}G(Y(s,W),s) \right) ds
\end{equation}
is also a martingale with respect to $\Im_t$ as soon as $G(Y,t)$ is twice continuously differentiable in $t$ and
$\sup_{Y,t} \left| \frac{\partial}{\partial t} G(Y,t) \right| + 
\left| \frac{\partial^2}{\partial t^2} G(Y,t) \right| < + \infty,$
which is condition (5.1) in \cite[App. 1 Sec. 5]{KipLan}.
\\ To proceed, consider $F$ of the form
\begin{equation}\label{FDEF}
 \Sigma_{\phi}(Y) := \sum_{f \in \MaxFacets(Y)} \phi(f)
\end{equation}
where, recall, $\MaxFacets(Y)$ are the maximal facets of $Y$ (the I-segments in the two-dimensional case)
whereas $\phi(\cdot)$ is a generic bounded and measurable functional on $(d-1)$-dimensional facets in $W,$
that is to say a bounded and measurable function on the space of closed $(d-1)$-dimensional polytopes in $W,$
possibly chopped off by the boundary of $W,$ with the standard measurable structure inherited from space of
closed sets in $W.$ Whereas the so-defined $F$ is not bounded and thus (\ref{MART}) cannot be applied
directly, we can apply it for $F_N := (F \wedge N) \vee -N,\; N \in {\Bbb N}$ which is bounded and
let $N\to\infty$
to conclude that $F(Y(t,W))-\int_0^t {\Bbb L}F(Y(s,W))ds$ is a local $\Im_t$-martingale, see Definition 5.15 in
\cite{KS} and take $T_N = \inf_{t\geq 0} |F(Y(t,W))| \geq N$ there. Now, apply the proof of Lemma 1
in \cite{NW05} where the number of cells in $Y(t,W),$ and hence for all $Y(s,W),\; s \leq t,$
is bounded by a Furry-Yule-type
linear birth process whose cardinality at any given finite time admits moments of all orders, to
conclude that $(F(Y(t,W))-\int_0^t {\Bbb L}F(Y(s,W))ds)_{t \leq a}$ is of class DL for all $a > 0$
in the sense of Definition 4.8 in \cite{KS}. Using now Problem 5.19 (i) in \cite{KS} we finally
conclude that $F(Y(t,W))-\int_0^t {\Bbb L}F(Y(s,W))ds$ is a martingale. 
Thus, applying (\ref{GEN}) and (\ref{MART}) for $F \equiv \Sigma_{\phi}$ we see that
\begin{equation}\label{EXPECT}
 \Sigma_{\phi}(Y(t,W)) - \int_0^t \int_{[W]} \sum_{x \in \iota(H;Y)} \phi(\Facet(x,H|Y(s,W))) \Lambda(dH)ds
\end{equation}
is a martingale with respect to $\Im_t$.

\section{First-Order Properties and Typical I-Faces}\label{secFIRST}
In this section we establish a number of first-order properties of $Y(t,W)=Y(t,\Lambda,W)$ for general locally finite non-atomic measure $\Lambda$, essentially obtained
 by comparison with suitable mixtures of Poisson hyperplane tessellations. Many of these properties
 are known in the translation-invariant set-up. Afterwards, these results are used to calculate several mean values for stationary and stationary and isotropic STIT tessellations $Y(t)$.
 
\subsection{Distribution of Cells and Lower Dimensional Faces}\label{seclowerdimfaces}
 
 The key to our results is formula (\ref{EXPECT})
 from Section \ref{secMARTINGALE}. To exploit it, consider the random measures 
 \begin{equation}\label{CellMeasure}
   {\cal M}^{Y(t,W)} := \sum_{c \in \Cells(Y(t,W))} \delta_c,\;\; {\Bbb M}^{Y(t,W)} := {\Bbb E}{\cal M}^{Y(t,W)}
 \end{equation}
 with $\delta_c$ standing for the unit Dirac mass at $c.$ In full analogy, define ${\cal M}^{\PHT(t\Lambda,W)}$ 
 and ${\Bbb M}^{\PHT(t\Lambda,W)}$ where $\PHT(t\Lambda,W)$ is the Poisson hyperplane tessellation
 with intensity measure $t\Lambda,$ restricted to $W.$  Further, put
 \begin{equation}\label{FK}
  {\cal F}_k^{Y(t,W)} := \sum_{f \in \MaxFaces_k(Y(t,W))} \delta_f,\;\; 
  {\Bbb F}_k^{Y(t,W)} := {\Bbb E}{\cal F}_k^{Y(t,W)},\; k=1,\ldots,d-1,
 \end{equation}
 where, recall, $\MaxFaces_k(Y)$ is the collection of $k$-dimensional I-faces of 
 $Y$ and, likewise, define
 $$  {\cal F}_k^{\PHT(t\Lambda,W)} := \sum_{f \in \Faces_k(\PHT(t\Lambda,W))} \delta_f,\;\; 
       {\Bbb F}_k^{\PHT(t\Lambda,W)} := {\Bbb E}{\cal F}_k^{\PHT(t\Lambda,W)},\; k=1,\ldots,d-1. $$
 Our first claim is
 \begin{theorem}\label{Meq}
  We have
  $$ {\Bbb M}^{Y(t,W)} = {\Bbb M}^{\PHT(t\Lambda,W)}. $$
 \end{theorem}
 Note, that in the particular case of $\Lambda$ being translation-invariant this reduces
 to the known fact that the typical cell of $Y(t)=Y(t\Lambda)$ coincides with the typical
 cell of $\PHT(t\Lambda).$ 
 \paragraph{Proof of Theorem \ref{Meq}}
 Using (\ref{GEN}) and (\ref{MART}) with 
 $$ F(Y) := \sum_{c \in \Cells(Y)} \phi(c) $$
 for general bounded measurable cell functional $\phi,$ with localization argument as
 the one preceding (\ref{EXPECT}) we conclude that
 $$  \int \phi d{\cal M}^{Y(t,W)} - 
     \int_0^t \int_{[W]} \sum_{x \in \iota(H;Y)} 
     [\phi(\Cell^+(x,H|Y(t,W))) + $$
\begin{equation}\label{Mmart1}
    \phi(\Cell^-(x,H|Y(t,W))) - \phi(\Cell(x,H|Y(t,W)))]
     \Lambda(dH) ds
 \end{equation}
 is a $\Im_t$-martingale.
 For a polyhedral cell $c \subseteq W,$ possibly chopped off by the boundary of $W,$
 and for $H \in [c]$ we write $c^+(H)$ and $c^-(H)$ to denote the
 cells into which $c$ gets divided by $H,$ lying respectively on the positive and negative
 side of $H.$ With this notation, (\ref{Mmart1}) says that
 $$ \int \phi d{\cal M}^{Y(t,W)} - 
     \int_0^t \int \int_{[c]} [\phi(c^+(H)) + \phi(c^-(H)) - \phi(c)] \Lambda(dH) {\cal M}^{Y(s,W)}(dc)
     ds $$
 is a $\Im_t$-martingale. Taking expectations leads to
 \begin{equation}\label{EQNforM}
  \int \phi d{\Bbb M}^{Y(t,W)} = \int_0^t \int \int_{[c]} [\phi(c^+(H)) + \phi(c^-(H)) - \phi(c)]
  \Lambda(dH) {\Bbb M}^{Y(s,W)}(dc) ds
 \end{equation}
 for all bounded measurable $\phi.$ To proceed, we regard ${\Bbb M}^{Y(s,W)}$ as an element of the space of bounded variation Borel
 measures on the family of polyhedral sub-cells of $W$ endowed with the standard measurable
 structure inherited from the space of closed sets in $W.$ Consider the linear
 operator $T_{\Lambda}$ on this measure space, given by
 \begin{equation}\label{TM}
  T_{\Lambda}(\mu) = \int \int_{[c]} [\delta_{c^+(H)} + \delta_{c^-(H)} - \delta_{c}] \Lambda(dH)
  \mu(dc).
 \end{equation}
 By the definition (\ref{TM}), $\left\|T_{\Lambda}(\mu)\right\|_{\rm TV} \leq (\int_{[W]} d\Lambda) \left\|\mu\right\|_{\rm TV}
 = \Lambda([W]) \left\|\mu\right\|_{\rm TV}$ where $\left\|\cdot\right\|_{\rm TV}$ is the standard total variation norm of a measure, see \cite[Def. 3.1.4]{BOG}. This
 inequality turns into equality when $\mu = \delta_{W}.$  Consequently, $T_{\Lambda}$ is a bounded
 operator of operator norm $\Lambda([W]) < +\infty.$
 The relation (\ref{EQNforM}) can be rewritten in differential form
 \begin{equation}\label{RNIEr}
  \frac{d}{dt} {\Bbb M}^{Y(t,W)} = T_{\Lambda} {\Bbb M}^{Y(t,W)},\;\; {\Bbb M}^{Y(0,W)} = \delta_{W}
 \end{equation}
 which, in view of the above properties of $T_{\Lambda},$ admits by standard theory (cf. \cite[IX.\S 2, Sec. 2]{Kato}) the unique
 solution 
 \begin{equation}\label{ROZW}
  {\Bbb M}^{Y(t,W)} = \exp(t T_{\Lambda}) \delta_{W},\;\; t \geq 0.
 \end{equation}   
 It is easily seen that exactly the same equations (\ref{EQNforM}),(\ref{RNIEr}) and thus also
 (\ref{ROZW}) hold for ${\Bbb M}^{\PHT(t\Lambda,W)}.$
 In particular, ${\Bbb M}^{Y(t,W)} = {\Bbb M}^{\PHT(t\Lambda,W)}$ as required.
 This completes the proof. $\hfill\Box$\\ \\ Having characterized ${\Bbb M}^{Y(t,W)}$ we now turn to ${\Bbb F}_k^{Y(t,W)}.$
 \begin{theorem}\label{FKeq}
  For all $k=1,\ldots,d-1$ we have
  $$ {\Bbb F}_k^{Y(t,W)} = (d-k) 2^{d-k-1} \int_0^t \frac{1}{s} {\Bbb F}_k^{\PHT(s\Lambda,W)} ds. $$
 \end{theorem}
  \paragraph{Proof of Theorem \ref{FKeq}}
   Fix $k \in \{1,\ldots,d-1\}.$ Let $\psi$ be a general bounded measurable function of a 
   $k$-dimensional I-face, as usual regarded as a closed subset of $W,$
   and for a $(d-1)$-dimensional I-facet $h$ put
   \begin{equation}\label{PhiPsi}
    \phi(h) := \sum_{f \in \Faces_k(h)} \psi(f) 
   \end{equation}
   noting that the $k$-dimensional I-faces of the tessellation $Y(t,W)$ are precisely the $k$-faces of its I-facets.  
   Using (\ref{EXPECT}), taking expectations and recalling (\ref{CellMeasure}) we see that
   $$ {\Bbb E} \Sigma_{\phi}(Y(t,W)) = \int \phi d{\Bbb F}_{d-1}^{Y(t,W)} =
        \int_0^t \int \int_{[c]} \phi(c \cap H) \Lambda(dH) {\Bbb M}^{Y(s,W)}(dc) ds. $$
  Applying Theorem \ref{Meq} we get
  $$ \int \phi d{\Bbb F}_{d-1}^{Y(t,W)} = \int_0^t \int \int_{[c]} \phi(c \cap H) 
      \Lambda(dH) {\Bbb M}^{\PHT(s\Lambda,W)}(dc) ds. $$
  However, applying Slivnyak's theory, see e.g. \cite[Thm 1.15]{SW2}, we obtain
  $$ \int \int_{[c]} \phi(c \cap H) \Lambda(dH) {\Bbb M}^{\PHT(\Lambda,W)}(dc)
      = \int \phi d{\Bbb F}_{d-1}^{\PHT(\Lambda,W)} $$
  and thus, upon taking $s\Lambda$ in place of $\Lambda,$ more generally,
  $$ \int \int_{[c]} \phi(c \cap H) \Lambda(dH) {\Bbb M}^{\PHT(s\Lambda,W)}(dc)
      = \frac{1}{s} \int \phi d{\Bbb F}_{d-1}^{\PHT(s\Lambda,W)} $$
  whence
  $$ \int \phi d{\Bbb F}_{d-1}^{Y(t,W)} = \int_0^t \frac{1}{s} \int \phi d{\Bbb F}_{d-1}^{\PHT(s\Lambda,W)} ds$$
  follows. Now note that, by (\ref{PhiPsi}),
  $$ \int \phi d{\Bbb F}_{d-1}^{Y(t,W)} = \int \psi d{\Bbb F}_k^{Y(t,W)} $$
 because each k-dimensional I-face is a k-face of precisely one I-facet in $Y(t,W).$ 
 Moreover, 
 $$ \int \phi d{\Bbb F}_{d-1}^{\PHT(s\Lambda,W)} = (d-k) 2^{d-k-1} \int \psi d{\Bbb F}_k^{\PHT(s\Lambda,W)} $$
 because each k-face of $\PHT(s\Lambda,W)$ is a k-face of $(d-k) 2^{d-k-1}$ facets of $\PHT(s\Lambda,W)$, see Theorems 10.1.2 and 10.3.1 in \cite{SW}.
 Hence,  we conclude that
 $$\int \psi d{\Bbb F}_{k}^{Y(t,W)} = (d-k) 2^{d-k-1} \int_0^t \frac{1}{s} \int \psi d{\Bbb F}_k^{\PHT(s\Lambda,W)} ds$$
for all $\psi$ bounded and measurable,  which completes the proof of the Theorem. $\hfill\Box$\\ \\
Some of our argument in the sequel will require a straightforward formal extension of Theorem \ref{FKeq}. Namely, we formally {\it mark} all I-facets of the tessellation $Y(t,W)$ by their {\it birth times}.
This gives rise to the birth-time augmented tessellation $\hat{Y}(t,W)$ with birth-time-marked I-facets
and makes the MNW-construction of $\hat{Y}(t,W)$ into a Markov process whose generator $\hat{{\Bbb L}}$
is a clear modification of ${\Bbb L}$ as given in (\ref{GEN}):
$$ \hat{\Bbb L} \hat{F}(\hat{Y}) = \int_{[W]} \sum_{x \in \iota(H;Y)} 
     [\hat{F}(\hat{Y} \cup [\Facet(x;H|Y),s]) - \hat{F}(\hat{Y})] \Lambda(dH) $$
for $\hat{F}$ bounded measurable on the space of birth time-marked tessellations of $W.$ Consequently, 
writing $\hat{\Bbb F}^{Y(t,W)}_k,\; k=1,\ldots,d-1,$ for the birth-time-marked version of
${\Bbb F}^{Y(t,W)}_k$ where each k-dimensional I-face is marked with its birth time, by
a straightforward modification of the proof of Theorem \ref{FKeq} we are led to
\begin{corollary}\label{TIMEMARKEDFQeq}
 For all $k=1,\ldots,d-1$ we have
 $$ \hat{\Bbb F}^{Y(t,W)}_k = (d-k) 2^{d-k-1} \int_0^t \frac{1}{s} 
      \left[ {\Bbb F}^{\PHT(s\Lambda,W)}_k \otimes \delta_s \right] ds. $$
\end{corollary} 

\subsection{Distribution of typical I-Faces and Related Mean Values}\label{sectypifaces}

\begin{figure}[t]
\begin{center}
 \includegraphics[width=12cm]{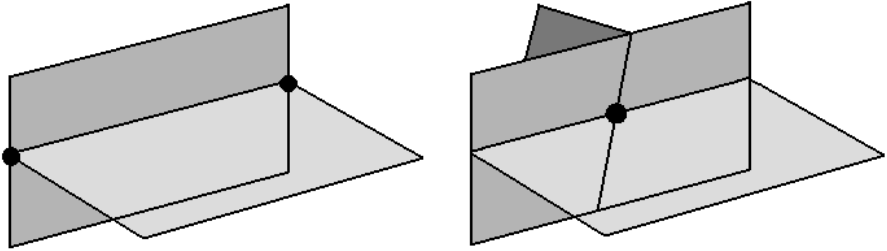}\\
 \caption{Vertices (left hand side) and crossing nodes (right hand side) that appear in $3$-dimensional STIT tessellations}\label{Fig2}
\end{center}
\end{figure}

We are now going to apply the results obtained in the last section to the stationary and isotropic set-up, i.e. with $\Lambda=\Lambda_{iso}$ and to calculate the distribution and the mean $k$-volume as well as other related mean values of the typical $k$-dimensional I-face $I_k^{(d)},\; k =0,\ldots,d-1,$ of the STIT tessellation $Y(t)$ in ${\Bbb R}^d$. In this paper we will restrict our attention only to proper $k$-dimensional I-faces, which we call \textit{$k$-dimensional I-faces} for short, see Section \ref{secSTIT} above for definitions. However, it should be noted that starting with dimension $3$, there appear also the so-called \textit{crossing $k$-faces} arising as intersections of traces of higher-dimensional I-faces inside another face. For $d=3$ these crossing k-faces appear only in the form of nodes. The $0$-dimensional I-faces are simply the vertices and the other type of 0-faces (nodes) are called crossing nodes. Figure \ref{Fig2} should clarify the situation. In higher dimensions the structure of crossing k-faces is even more complicated, with
the notable exception of the case $k=d-1$ where there only exists one type of I-faces, namely the I-facets.
As stated above, throughout this paper whenever talking about $k$-dimensional I-faces we always mean the proper ones.
To proceed, let $\varphi_k:\MaxFaces_k\rightarrow{\Bbb R}$ be a translation-invariant, non-negative measurable function for $1\leq k\leq d-1$ and denote by $\overline{\varphi}_k(Y(t))$ the (possibly infinite) $\varphi_k$-density of $Y(t)$ in the sense of \cite[Chap. 4.1]{SW}, i.e. $$\overline{\varphi}_k(Y(t))=\lim_{r\rightarrow\infty}{1\over r^d\Vol_d(W)}{\Bbb E}\sum_{f\in\MaxFaces_k(Y(t,rW))}\varphi_k(f)$$ with $W\subset{\Bbb R}^d$ some bounded convex set with positive and finite volume. The existence of this limit is guaranteed by Thm 4.1.3 ibidem. Using now Theorem \ref{FKeq} from above and Campbell's formula, we obtain, possibly with both sides infinite,
\begin{eqnarray}
\nonumber \overline{\varphi}_k(Y(t)) &=& \lim_{r\rightarrow\infty}{1\over r^d\Vol_d(W)}{\Bbb E}\sum_{f\in\MaxFaces_k(Y(t,rW))}\varphi_k(f)\\
\nonumber &=& \lim_{r\rightarrow\infty}{1\over r^d\Vol_d(W)}{\Bbb E}\int\varphi_k(f){\cal F}_k^{Y(t,rW)}(df)\\
\nonumber &=& \lim_{r\rightarrow\infty}{1\over r^d\Vol_d(W)}\int\varphi_k(f){\Bbb F}_k^{Y(t,rW)}(df)\\
\nonumber &=& \lim_{r\rightarrow\infty}{1\over r^d\Vol_d(W)}(d-k)2^{d-k-1}\int_0^t{1\over s}\left[\int\varphi_k(f){\Bbb F}_k^{\PHT(s,rW)}(df)\right]ds\\
\nonumber &=& (d-k)2^{d-k-1}\int_0^t{1\over s}\left[\lim_{r\rightarrow\infty}{1\over r^d\Vol_d(W)}{\Bbb E}\sum_{f\in\Faces_k(\PHT(s,rW))}\varphi_k(f)\right]ds\\
&=& (d-k)2^{d-k-1}\int_0^t{1\over s}\overline{\varphi}_k(\PHT(s))ds,\label{EQMIXFK}
\end{eqnarray}
where $\PHT(s)$ ($\PHT(s,rW)$) is the stationary and isotropic Poisson hyperplane tessellation with surface intensity $s$ (restricted to the window $rW$). We use the formula (\ref{EQMIXFK}) with $\varphi_k\equiv 1$ and denote in this case $\overline{\varphi}_k(Y(t))=:N_{k,I}^{(d)}$, referred to as the
\textit{intensity of $k$-dimensional I-faces}. Then we have upon applying \cite[Thm. 10.3.3]{SW} together with (\ref{EQMIXFK}) 
\begin{equation}\label{NKIEXPR}
N_{k,I}^{(d)}=(d-k)2^{d-k-1}\int_0^t{1\over s}{d\choose k}\kappa_d\left({\kappa_{d-1}\over d\kappa_d}\right)^ds^dds=(d-k)2^{d-k-1}{\kappa_d\over d}{d\choose k}\left({\kappa_{d-1}\over d\kappa_d}\right)^dt^d.
\end{equation}
For a $k$-dimensional polytope $f$, let $c(f)$ be some associated center function with the property that $c(f+x)=c(f)+x$ for any $x\in{\Bbb R}^d$ (take for example the Steiner point or the center of the minimal circumscribed ball of $f$). Now, we use again (\ref{EQMIXFK}) with $\varphi_k(f):={\bf 1}[\cdot](f-c(f))$ together with \cite[Eq. (4.8,4.9)]{SW}
and obtain in this case 
\begin{equation}\label{MIXTYPEQ}
N_{k,I}^{(d)}{\Bbb Q}_k^{Y(t)}=(d-k)2^{d-k-1}\int_0^t{1\over s}\gamma_{\PHT(s)}^{(k)}{\Bbb Q}_k^{\PHT(s)}ds,
\end{equation}
where $\gamma_{\PHT(s)}^{(k)}$ is the mean number of $k$-faces of $\PHT(s)$ per unit volume, ${\Bbb Q}_k^{\PHT(s)}$ is the distribution of the typical $k$-face of $\PHT(s)$ and ${\Bbb Q}_k^{Y(t)}$ is the distribution of the typical $k$-dimensional I-face of the STIT tessellation $Y(t)$, see \cite[Thm. 4.1.1]{SW} and the discussion ibidem for definition of typical faces and typical grains of particle processes in general. Inserting now the value for $N_{k,I}^{(d)}$ and the one for $\gamma_{\PHT(s)}^{(k)}$ from \cite[Thm. 10.3.3]{SW} and noting that $N_{k,I}^{(d)}={d-k\over d}2^{d-k-1}\gamma_{\PHT(t)}^{(k)}$ by comparing the total masses in both sides of (\ref{MIXTYPEQ}), we arrive at
\begin{theorem}\label{CORDISTR} The distribution ${\Bbb Q}_k^{Y(t)}$ of the typical $k$-dimensional I-face of $Y(t)$ is given by $${\Bbb Q}_k^{Y(t)}=\int_0^t{ds^{d-1}\over t^d}{\Bbb Q}_k^{\PHT(s)}ds,$$ where ${\Bbb Q}_k^{\PHT(s)}$ is the distribution of the typical k-face of the stationary and isotropic Poisson hyperplane tessellation $\PHT(s)$ with surface intensity $s$.
\end{theorem}
The preceding theorem can be rephrased by saying that the distribution of the \textit{typical $k$-dimensional I-face} of a STIT tessellation is a mixture of suitable rescalings of the distribution of the typical $k$-dimensional face of Poisson hyperplane tessellations.\\ In what follows, $I_k^{(d)}$ stands for the $k$-polytope with distribution ${\Bbb Q}_k^{Y(t)}$. For the typical I-segment $I_1^{(d)}$ of the stationary and isotropic STIT tessellation $Y(t)$ in ${\Bbb R}^d$ we can obtain even more than the statement of Theorem \ref{CORDISTR}, namely the fact that the length distribution of $I_1^{(d)}$ is a mixture of exponential distributions, answering thereby a question formulated in the Outlook of \cite{M09}. It is also interesting to note that formally marking the $k$-dimensional I-faces with
their birth-times and repeating the argument leading to Theorem \ref{CORDISTR} with Theorem \ref{FKeq}
replaced by its time-marked extension in Corollary \ref{TIMEMARKEDFQeq} we obtain the birth time-marked
extension of Theorem \ref{CORDISTR}
\begin{corollary}\label{TIMEMARKEDCORDISTR}
 The distribution $\hat{\Bbb Q}_k^{Y(t)}$ of the typical birth-time-marked $k$-dimensional I-face of $Y(t)$ is given by 
 $$\hat{\Bbb Q}_k^{Y(t)}=\int_0^t{ds^{d-1}\over t^d}\left[ {\Bbb Q}_k^{\PHT(s)} \otimes \delta_s \right] ds.$$
\end{corollary} 
To proceed, using Theorem \ref{CORDISTR}, 
for the \textit{length density} $p_l^{(d)}(x)$ of the typical I-segment of a stationary and isotropic STIT tessellation $Y(t)$ in ${\Bbb R}^d$ we obtain because of key property (b) from Section \ref{secSTIT}
\begin{equation} p_l^{(d)}(x) = \int_0^t\lambda_1se^{-\lambda_1sx}{ds^{d-1}\over t^d}ds={d\over(\lambda_1t)^dx^{d+1}}\gamma(d+1,\lambda_1tx),\label{eqPLDI1D}\end{equation}
where $\gamma(\cdot,\cdot)$ is the lower incomplete Gamma-function and $\lambda_1$ is given by (\ref{LAMBDAK}). From this it is easily seen that for the length of the typical I-segment only the moments of order $1$ to $d-1$ exist, see Section \ref{secHIGERMOMENTS} below. In particular for $d=2$ and $d=3$ we have the densities
\begin{eqnarray}
p_l^{(2)}(x) &= &{1\over t^2x^3}\left(\pi^2-(\pi^2+2\pi tx+2t^2x^2)e^{-{2\over\pi}tx}\right),\label{EQPLANARPL}\\
\nonumber p_l^{(3)}(x) &=& {3\over t^3x^4}\left(48-(48+24tx+6t^2x^2+t^3x^3)e^{-{1\over 2}tx}\right).
\end{eqnarray}
We use now once more (\ref{EQMIXFK}) but this time with $\varphi_k(f)=V_j(f)$ being the $j$-th intrinsic volume of the $k$-polytope $f$, $0\leq j\leq k.$ Define $\overline{\varphi}_k(Y(t))=:d_j^{(k,d)}$ and use \cite[Thm 10.3.3]{SW} to obtain $$d_j^{(k,d)}=2^{d-k-1}{d-k\over d-j}{d-j\choose d-k}{d\choose j}\left({\kappa_{d-1}\over d\kappa_d}\right)^{d-j}{\kappa_d\over\kappa_j}t^{d-j}$$ as the value for the \textit{density of the $j$-th intrinsic volume of the (proper) $k$-skeleton} of the stationary and isotropic STIT tessellation $Y(t)$. For $j=k$ this yields 
\begin{equation}\label{SVEXPR}
 S_V^{(k,d)}=d_k^{(k,d)}=2^{d-k-1}{d\choose k}{\kappa_d\over\kappa_k}\left({\kappa_{d-1}\over d\kappa_d}\right)^{d-k}t^{d-k},
\end{equation}
 i.e. $S_V^{(k,d)}$ is the mean $k$-volume of (proper) $k$-dimensional I-faces of $Y(t)$ per unit volume. Note, that in particular $$ N_0^{(d)} := N_{0,I}^{(d)}=S_V^{(0,d)}=2^{d-1}{\kappa_d}\left({\kappa_{d-1}\over d\kappa_d}\right)^dt^d$$ is the mean number of vertices of $Y(t)$ per unit volume in agreement with (\ref{NKIEXPR}) above. Using the identity $$N_{k,I}^{(d)}\cdot{\Bbb E}\Vol_k(I_k^{(d)})=S_V^{(k,d)}$$ combined with (\ref{NKIEXPR}) and (\ref{SVEXPR}) leads to $${\Bbb E}\Vol_k(I_k^{(d)})={S_V^{(k,d)}\over N_{k,I}^{(d)}}={d\over(d-k)\kappa_k}\left({d\kappa_d\over\kappa_{d-1}}\right)^k{1\over t^k}$$ for the mean $k$-volume of the typical $k$-dimensional I-face $I_k^{(d)}$. Especially for $k=1$ and $k=d-1$ we have $${\Bbb E}\Vol_1(I_1^{(d)})={d\sqrt{\pi}\Gamma\left({d-1\over 2}\right)\over 2\Gamma\left({d\over 2}\right)}{1\over t},\ \ \ {\Bbb E}\Vol_{d-1}(I_{d-1}^{(d)})={d\over\kappa_{d-1}}\left({2\sqrt{\pi}\Gamma\left({d+1\over 2}\right)\over\Gamma\left({d\over 2}\right)}\right)^{d-1}{1\over t^{d-1}}.$$ In the planar case our formulas specialize to $N_0^{(2)}={2\over\pi}t^2$, $N_{1,I}^{(2)}={1\over\pi}t^2$ and $L_I={\Bbb E}\Vol_1(I_1^{(2)})={\pi\over t}$ and for space dimension $d=3$ we reproduce the results from \cite{TW}, namely $N_0^{(2)}={\pi\over 12}t^3$ (this is the mean number of T-vertices in the sense of \cite{TW} per unit volume), $N_{1,I}^{(2)}={\pi\over 12}t^3$, $N_{2,I}^{(3)}={\pi\over 48}t^3$ and $L_V=S_V^{(1,3)}={\pi\over 4}t^2$. Thus, the mean length of the typical I-segment in the $3$-dimensional case equals ${\Bbb E}\Vol_1(I_1^{(3)})={3\over t}$ and the mean area of the typical I-face is given by ${\Bbb E}\Vol_2(I_2^{(3)})={48\over\pi t^2}$.\\ Theorem \ref{CORDISTR} may also be used to calculate the \textit{mean f-vector} $(f^{(k)}_0,f^{(k)}_1,\ldots,f^{(k)}_{k-1})$ of the typical $k$-dimensional I-face $I_k^{(d)}$, i.e. $f^{(k)}_j$ is the mean number of $j$-dimensional I-faces of the typical k-dimensional I-face $I_k^{(d)}$ (thus $f^{(k)}_0$ is the mean number of vertices of a $k$-face) for $0\leq j<k\leq d$, where for $k=d$ we abuse the notation and let $I_d^{(d)}$ stand for the typical cell of $Y(t)$. Theorem \ref{CORDISTR} together with the crucial property (a) of STITS tessellations immediately imply that this vector is the same as the mean f-vector of a $k$-dimensional Poisson polytope. From \cite[Thm. 10.3.1]{SW} we deduce now $$f^{(k)}_j=2^{k-j}{k\choose j}.$$ In particular, the typical cell of $Y(t)$ has in the mean $2^d$ vertices, whereas the zero cell has ${d!\kappa_d^2\over 2^d}$, cf. \cite[Thm. 10.4.9]{SW}. For the typical cell in the case $d=2$ this yields the mean f-vector $(4,4)$ and in the spatial case $d=3$ we obtain for the typical $2$-dimensional I-facet and the typical cell the mean f-vectors $(4,4)$ and $(8,12,6)$, respectively, which is already known from \cite{TW}. We have thus proved
\begin{corollary}\label{corfvector} The mean $k$-volume of the typical $k$-dimensional I-face is given by \begin{equation}{\Bbb}\Vol_k(I_k^{(d)})={d\over(d-k)\kappa_k}\left({d\kappa_d\over\kappa_{d-1}}\right)^k{1\over t^k}={d\over(d-k)\kappa_k}\left({2\sqrt{\pi}\Gamma\left({d+1\over 2}\right)\over\Gamma\left({d\over 2}\right)}\right)^k{1\over t^k}\label{eqMEANVIKD}\end{equation} and the mean f-vector of $I_k^{(d)}$ equals $(f^{(k)}_0,\ldots,f^{(k)}_{k-1})$ with $f^{(k)}_j=2^{k-j}{k\choose j}$ for $0\leq j\leq k-1<d$.
\end{corollary}
Beside the typical geometry of the boundary of the typical $k$-dimensional I-face, some information may also be deduced about the geometry of the relative interior of $I_k^{(d)}$. Denote by  $N_{I,k,j}^{(d)}$ the the mean number of (proper) $j$-dimensional I-faces in the relative interior of the typical $k$-dimensional I-face $I_k^{(d)}$ and observe that this mean value is given by \begin{equation} N_{I,k,j}^{(d)}={N_{j,I}^{(d)}\over N_{k,I}^{(d)}}={(d-j){d\choose j}\over(d-k){d\choose k}}2^{k-j},\ \ \ \ 0\leq j<k<d,\label{NIKJ}\end{equation}
because each such j-dimensional I-face is contained in the relative interior of precisely one k-dimensional I-face.
 For $k=d-1$, i.e. when regarding I-facets, the formula specializes to $$N_{I,d-1,j}^{(d)}={d-j\over d}{d\choose j}2^{d-1-j}={d-1\choose j}2^{d-1-j}$$ and we obtain the mean value formula \begin{equation}N_{I,d-1,j}^{(d)}=f_j^{(d-1)}.\label{MVF2}\end{equation} For example we have in the planar case $N_{I,1,0}^{(d)}=2$ in accordance with \cite{NW06} and for $d=3$, $k=2$ and $j=1$ the value $N_{I,2,1}^{(3)}=4$, which fits with the results obtained in \cite{TW}.\\ Our next goal is a formula for the \textit{specific (or mean) $j$-th intrinsic volume} ${\Bbb E}V_j(I_k^{(d)})$, $0\leq j\leq k$, in the sense of \cite{SW} of the typical $k$-dimensional I-face $I_k^{(d)}$. In particular ${\Bbb E}V_k(I_k^{(d)})$ will again be the mean $k$-volume, $2{\Bbb E}V_{k-1}(I_k^{(d)})$ the mean surface area, ${2\kappa_{k-1}\over k\kappa_k}{\Bbb E}V_1(I_k^{(d)})$ the mean breath and, trivially,  ${\Bbb E}V_0(I_k^{(d)}) \equiv 1$ the mean Euler-number of the typical $k$-dimensional I-face $I_k^{(d)}$. Such a formula for the specific intrinsic volumes could easily be obtained by dividing their intensities $d_j^{(k,d)}$ by the face intensities $N_{k,I}^{(d)}$. For explanatory reasons we will give another proof though, highlighting the temporal construction of STIT tessellations and the
underlying ideas.
\begin{corollary}\label{corMEANINTVOL} The specific $j$-th intrinsic volume of the typical $k$-dimensional I-face of a stationary and isotropic STIT tessellation with surface intensity $t$ is given by \begin{equation}{\Bbb E}V_j(I_k^{(d)})={d\over (d-j)\kappa_j}{k\choose j}\left({2\sqrt{\pi}\Gamma\left({d+1\over 2}\right)\over\Gamma\left({d\over 2}\right)}\right)^j{1\over t^j}.\label{eqEVJ}\end{equation}
\end{corollary}
\paragraph{Proof of Corollary \ref{corMEANINTVOL}}
In analogy to what we did in Corollaries \ref{TIMEMARKEDFQeq} and \ref{TIMEMARKEDCORDISTR},
we regard the collection of all $k$-dimensional I-faces of $Y(t)$ as a marked point process of $k$-faces in the space of $k$-polytopes, where the marks are given by the birth times of the faces and the birth time of a $k$-face is defined to be the birth time of the unique $(d-1)$-facet which induces the birth of the $k$-face. Given now the birth time $\beta=s$ of $I_k^{(d)},$
the conditional specific $j$-th intrinsic volume of the typical $k$-dimensional I-face $I_k^{(d)}$ is  $${\Bbb E}[V_j(I_k^{(d)})|\beta=s]={k\choose j}\left({k\kappa_k\over\kappa_{k-1}}\right)^j{1\over\kappa_j(\lambda_k s)^j}={1\over\kappa_j}{k\choose j}\left({2\sqrt{\pi}\Gamma\left({d+1\over 2}\right)\over\Gamma\left({d\over 2}\right)}\right)^j{1\over s^j}$$ with $\lambda_k$ as in (\ref{LAMBDAK}). Here we have used Crofton's formula for random tessellations \cite[Eq. (10.28)]{SW} and key properties (a) and (b) for stationary and isotropic STIT tessellations. Now, Theorem \ref{CORDISTR} and Corollary \ref{TIMEMARKEDCORDISTR} say in this context that the the birth time density of $I_k^{(d)}$ is given by $${ds^{d-1}\over t^d},\ \ \ 0<s<t.$$ Integration of ${\Bbb E}[V_j(I_k^{(d)})|\beta=s]$ with respect to this birth time density yields the value for ${\Bbb E}V_j(I_k^{(d)})$, i.e.
\begin{eqnarray}
\nonumber {\Bbb E}V_j(I_k^{(d)}) &=& \int_0^t{ds^{d-1}\over t^d}{\Bbb E}[V_k(I_k)|\beta=s]ds=\int_0^t{ds^{d-1}\over t^d}{1\over\kappa_j}{k\choose j}\left({2\sqrt{\pi}\Gamma\left({d+1\over 2}\right)\over\Gamma\left({d\over 2}\right)}\right)^j{1\over s^j}ds\\
\nonumber &=& {d\over (d-j)\kappa_j}{k\choose j}\left({2\sqrt{\pi}\Gamma\left({d+1\over 2}\right)\over\Gamma\left({d\over 2}\right)}\right)^j{1\over t^j}
\end{eqnarray}
and completes the proof. $\hfill\Box$\\ \\
For $j=k$, (\ref{eqEVJ}) specializes to formula (\ref{eqMEANVIKD}), for $j=0$ we have the constant value $1$, which is clear since $I_k^{(d)}$ is a convex set with probability $1$, and for $j=1$ and $j=k-1$ we obtain
\begin{eqnarray}
\nonumber {\Bbb E}V_1(I_k^{(d)}) &=& {kd\over d-1}{\sqrt{\pi}\Gamma\left({d+1\over 2}\right)\over\Gamma\left({d\over 2}\right)}{1\over t},\\
\nonumber {\Bbb E}V_{k-1}(I_k^{(d)}) &=& {kd\over (d-k+1)\kappa_{k-1}}\left({2\sqrt{\pi}\Gamma\left({d+1\over 2}\right)\over\Gamma\left({d\over 2}\right)}\right)^{k-1}{1\over t^{k-1}}.
\end{eqnarray}
It is our next goal to generalize some mean value formulas known from earlier papers for lower-dimensional stationary STIT tessellations to arbitrary space dimensions. For this reason we introduce the notion of \textit{$k$-dimensional J-faces}. By such a face for $Y(t)$ we mean {\it any} $k$-face of a $d$-dimensional cell, not necessarily maximal. Moreover, a J-face is counted with multiplicity given by the number of cells in whose boundary it lies. This way, the collection of $k$-dimensional J-faces with such multiplicities coincides with the collection of all k-faces of 
the tessellation cells in $Y(t),$ again with respective multiplicities. This convention, as adopted in the existing literature, has the clear and important advantage of ensuring that the typical $k$-dimensional J-face of $Y(t),$
denoted by  $J_k^{(d)}$ below, has the same distribution as the typical k-face of a corresponding Poisson hyperplane
tessellation $\PHT(t\Lambda)$ with the same intensity measure $t\Lambda,$ because of the Poisson typical cell
property (a) of STIT tessellations. Consequently,
\begin{equation}\label{TypicalJKFace}
  {\Bbb E} V_j(J_k^{(d)}) = {\Bbb E}V_j(\mbox{Typical k-Face}(\PHT(t\Lambda))),\; 0 \leq j \leq k < d. 
\end{equation}
A look at the display (\ref{eqEVJ}) in Corollary \ref{corMEANINTVOL} and a comparison with the
second formula on page 490 of \cite{SW} leads now from (\ref{TypicalJKFace}) to the following
comparison relation
\begin{equation} {\Bbb E}V_j(I_k^{(d)})={d\over d-j}{\Bbb E}V_j(J_k^{(d)}),\ \ \ 0\leq j\leq k<d,\label{MVF23}\end{equation} and especially for $j=k$ we have
 \begin{equation}\nonumber{\Bbb E}\Vol_k(I_k^{(d)})={d\over d-k}{\Bbb E}\Vol_k(J_k^{(d)}).\end{equation}
Writing now $N_{k,J}^{(d)}$ for the mean number of $k$-dimensional J-faces per unit volume we get  therefore
$$ {\Bbb E}\Vol_k(J_k^{(d)})\cdot N_{k,J}^{(d)}={d-k\over d}{\Bbb E}\Vol_k(I_k^{(d)})\cdot N_{k,J}^{(d)}. $$
On the other hand, taking into account that each point of the k-skeleton of $Y(t)$ belongs to precisely $(d-k+1)$ of
the k-dimensional J-faces we see that
$$ {\Bbb E}\Vol_k(J_k^{(d)})\cdot N_{k,J}^{(d)} = (d-k+1)S_V^{(k,d)}. $$
Putting these together we are led to 
\begin{eqnarray}
N_{k,J}^{(d)}&=&{d(d-k+1)\over d-k}N_{k,I}^{(d)}\label{MVF24}\\
\nonumber &=&(d-k+1)2^{d-k-1}\kappa_d{d\choose k}\left({\kappa_{d-1}\over d\kappa_d}\right)^dt^d.
\end{eqnarray}
For $d=2$, $k=1$ and $d=3$, $k=1$ or $k=2$ these relationships are known from \cite{NW06} and \cite{TW}. The value for $N_{k,J}^{(d)}$ can for example be used to obtain a formula for the mean number of vertices $N_{J,k,0}^{(d)}$ in the relative interior of the typical $k$-dimensional J-face: \begin{equation}N_{J,k,0}^{(d)}={N_0^{(d)}\over N_{k,J}^{(d)}}={2^k\over(d-k+1){d\choose k}}.\label{MVF3}\end{equation} This yields in the planar case $1\over 2$ and for $d=3$ and $k=2$ the value $2\over 3$, both known from \cite{NW06} and \cite{TW}, respectively. Also observe that the mean f-vector of $J_k^{(d)}$ is the same as the corresponding one for the typical I-face $I_k^{(d)}$.
\begin{remark} We wish to remark that the approach used in the proof of Corollary \ref{corMEANINTVOL} and a similar approach used in earlier papers on STIT tessellations, which is essentially based on marked point processes, is now fully justified by the considerations of this subsection, especially Corollary \ref{TIMEMARKEDCORDISTR}. The potential underlying \textit{Slivnyak-type theorem for I-faces and I-facets} is replaced here by our martingale technique and the fact that (\ref{MART}) is a martingale can in some sense be seen as such a Slivnyak-type result for STIT tessellations or more generally for iteration infinitely divisible random tessellations.
\end{remark}

\subsection{Higher Moments for Typical I-Faces}\label{secHIGERMOMENTS}

This section links mainly to the results from the last subsection and it is our goal to establish formulas for the second and third moment of the $k$-volume of the typical $k$-dimensional I-face of stationary and isotropic STIT tessellations $Y(t)$ with surface intensity $t$. Moreover, we will establish a formula for the second moment of the $(k-1)$-dimensional surface area of the boundary of $I_k^{(d)}$. For the special case $k=1$, formulas for the second or higher moments of the length of the typical I-segment can be obtained directly from the explicit density (\ref{eqPLDI1D}). For $1\leq k\leq d-1$ we calculate
\begin{eqnarray}
\nonumber {\Bbb E}\Vol_1^n(I_1^{(d)}) &=& \int_0^\infty x^n\int_0^t\lambda_1 se^{-\lambda_1 sx}{ds^{d-1}\over t^d}dsdx=\int_0^td\lambda_1{s^d\over t^d}\int_0^\infty x^ne^{-\lambda_1 sx}dxds\\
 &=& {d\cdot n!\over t^d\lambda_1^n}\int_0^t{s^d\over s^{n+1}}ds={d\cdot n!\over d-n}\left({\sqrt{\pi}\Gamma\left({d+1\over 2}\right)\over\Gamma\left({d\over 2}\right)}\right)^n{1\over t^n}.\label{eqnmomi1d}
\end{eqnarray}
For $k>1$, i.e. for higher dimensional I-faces, we start with the following general statement, which directly follows from Theorem \ref{CORDISTR} and the fact that $\TypicalCell(\PHT_k(\lambda_ks))$ has the same distribution as the typical $k$-face $\TKF(\PHT(s))$ of a stationary and isotropic Poisson hyperplane tessellation with surface intensity $s$ in ${\Bbb R}^d$, where by $\TypicalCell(\PHT_k(\lambda_ks))$ we mean the typical cell of a stationary and isotropic Poisson hyperplane tessellation in ${\Bbb R}^k$ with surface intensity $\lambda_ks$ and $\lambda_k$ given by (\ref{LAMBDAK}).
\begin{corollary}\label{cormomintvol} The $n$-th moment of the $j$-th intrinsic volume of the typical $k$-dimensional I-face $I_k^{(d)}$ of a stationary and isotropic STIT tessellation $Y(t)$ in ${\Bbb R}^d$ with surface intensity $t$ satisfies 
\begin{eqnarray}
\nonumber {\Bbb E}V_j^n(I_k^{(d)}) &=& \int_0^t{ds^{d-1}\over t^d}{\Bbb E}V_j^{n}(\TKF(\PHT(s)))ds\\
&=& \int_0^t{ds^{d-1}\over t^d}{\Bbb E}V_j^{n}(\TypicalCell(\PHT_k(\lambda_ks)))ds \label{eqSMC1}
\end{eqnarray}
with $0\leq j\leq k<d-1$. This value is finite if and only if $d-jn>0$.
\end{corollary}
At first, we apply (\ref{eqSMC1}) for $j=k $ and $n=2$. It is well known that in this case $${\Bbb E}\Vol_k^2(\TypicalCell(\PHT_k(\lambda_ks)))={k!\over 2^k}\left({k\kappa_k\over\kappa_{k-1}}\right)^{2k}{1\over (\lambda_ks)^{2k}}={k!\over 2^k}\left({d\kappa_d\over\kappa_{d-1}}\right)^{2k}{1\over s^{2k}},$$ see \cite[p. 179]{Mat} and we obtain by integration $${\Bbb E}\Vol_k^2(I_k^{(d)})={d\over d-2k}{k!\over 2^k}\left({2\sqrt{\pi}\Gamma\left({d+1\over 2}\right)\over\Gamma\left({d\over 2}\right)}\right)^{2k}{1\over t^{2k}},$$ whenever $d-2k>0$ and $+\infty$ in the other cases. For $d\geq 3$ and $k=1$ our formula yields the second moment of the length of the typical I-segment, see (\ref{eqnmomi1d}), which specializes for $d=3$ to the value ${\Bbb E}\Vol_1^2(I_1^{(3)})={24\over t^2}$.\\ In view of the results from \cite[Sec. 6.3]{Mat}, a similar formula is also available for the third moment of $\Vol_k(I_k^{(d)})$, i.e. for the case $j=k$ and $n=3$. For $d-3k>0$ we have $${\Bbb E}\Vol_k^3(I_k^{(d)})={d\over d-3k}2^{2k}\pi^{k-3\over 2}\left({2\sqrt{\pi}\Gamma\left({d+1\over 2}\right)\over\Gamma\left({d\over 2}\right)}\right)^{3k}{\Gamma\left(1+{k\over 2}\right)\Gamma\left(k+{3\over 2}\right)\Gamma\left({k+1\over 2}\right)^3\over\Gamma\left({3(k+1)\over 2}\right)}{1\over t^{3k}}$$ and ${\Bbb E}\Vol_k^3(I_k^{(d)})=+\infty$, whenever $d-3k\leq 0$. In particular for $d\geq 4$ we obtain once more the third moment of the length of the typical I-segment (\ref{eqnmomi1d}) with $n=3$. By the same method we also get for $d\geq 3$ and $k\geq 2$ a formula for the second moment of the surface area of the boundary $\bd(I_k^{(d)})$ of $I_k^{(d)}$ (this corresponds up to a factor $2$ to the case $j=k-1$ and $n=2$ in Corollary \ref{cormomintvol}):
\begin{eqnarray}
\nonumber {\Bbb E}\Vol_{k-1}^2(\bd(I_k^{(d)})) &=& {d\over d+2-2k}{k!\over 2^{k-2}}\left({d\kappa_d\over\kappa_{d-1}}\right)^{2-2k}\left(1+\left({k\kappa_k\over 2\kappa_{k-1}}\right)^2\right){1\over t^{2k-2}}\\
\nonumber &=& {d\over d+2-2k}{k!\over 2^{k-2}}\left({2\sqrt{\pi}\Gamma\left({d+1\over 2}\right)\over\Gamma\left({d\over 2}\right)}\right)^{2-2k}\left(1+\pi{\Gamma\left({k+1\over 2}\right)^2\over\Gamma\left({k\over 2}\right)^2}\right){1\over t^{2k-2}}.
\end{eqnarray}
if $d+2-2k>0$, see again \cite{Mat}. Especially for $d=3$ and $k=2$ we obtain the second moment of the perimeter of the typical $I$-facet $${\Bbb E}\Vol_1^2(\bd(I_2^{(3)}))={3\over 4}\left(1+{\pi^2\over 4}\right){1\over t^2}.$$ For the particular case $d=3$ and $k=2$ we also know that the second moment of the number of vertices of $I_2^{(3)}$ equals $${\pi^2\over 2}+12.$$ Explicit formulas for other or higher moments of the intrinsic volumes of $I_k^{(d)}$ are not possible, since the corresponding values for Poisson hyperplane tessellations are currently not known up to our best knowledge.

\subsection{Typical Geometry beyond the Isotropic Regime}\label{SECBEHINDISO}

It is our aim to explore in this section some aspects of the typical geometry of STIT tessellations and more generally iteration infinitely divisible random tessellations beyond the isotropic regime considered in the last paragraphs. We first put ourself in a very general position by taking $\Lambda$ to be an arbitrary non-atomic and locally finite measure on the space $\cal H$ of hyperplanes, which is absolutely continuous with respect to the motion-invariant measure $\Lambda_{iso}$. In this case the resulting tessellation $Y(t),$ regarded as a particle
process of convex cells in the sense of \cite[Chap. 4]{SW}, admits a translation-regular intensity measure in the sense
of \cite[Chap. 11]{SW} and the theory presented there applies. Furthermore, the intensity function $N_{k,I}^{(d)}(z),\; z \in {\Bbb R}^d$ of $k$-dimensional I-faces, the intensity function $d_j^{(k)}(z),\; z \in {\Bbb R}^d$ of the $j$-th intrinsic volume of the $k$-skeleton and the lower-dimensional directional distributions in the sense of \cite[Chap. 11]{SW} of the iteration infinitely divisible random tessellation $Y(t\Lambda)$ are well defined. Combining the results from \cite[Chap. 4 and Chap. 11]{SW} with Theorem \ref{FKeq} from above yields at first
\begin{corollary}\label{corDD} The directional distribution of the random process of $k$-dimensional I-faces of an iteration infinitely divisible random tessellation with translation-regular intensity measure $\Lambda$ coincides with the directional distribution of the $(d-k)$-th intersection process of a Poisson hyperplane tessellation with the same intensity measure.
\end{corollary}
In the translation-invariant case and only for $d=2$ this was already shown in \cite{M09} by a quite different argument.\\ Denote by $\gamma_{\PHT(s\Lambda)}^{(k)}(z)$ the intensity function of $k$-faces of a Poisson hyperplane tessellation with translation-regular intensity measure $s\Lambda$. Then we obtain by applying Theorem \ref{FKeq} above and \cite[Thm. 11.4.2]{SW} the \textit{local mean value formula}
\begin{eqnarray}
\nonumber N_{k,I}^{(d)}(z) &=& (d-k)2^{d-k-1}\int_0^t{1\over s}\gamma_{\PHT(s\Lambda)}^{(k)}(z)ds=(d-k)2^{d-k-1}{d\choose k}\int_0^t{1\over s}\gamma_{\PHT(s\Lambda)}^{(0)}(z)ds\\
\nonumber &=& {d-k\over d}2^{-k}{d\choose k}N_0^{(d)}(z)=2^{-k}{d-1\choose k}N_0^{(d)}(z)
\end{eqnarray}
for Lebesgue almost all $z\in{\Bbb R}^d$. Analogously, for the intensity function $d_j^{(k)}(z)$ of the $j$-th intrinsic volume of the $k$-skeleton of $Y(t\Lambda)$ we obtain the \textit{local mean value relation} $$d_j^{(k)}(z)={d-k\over d-j}2^{j-k}{d-j\choose d-k}d_j^{(j)}(z)$$ for almost all $z\in{\Bbb R}^d$ wrt. Lebesgue measure. In the translation-invariant set-up, the considered intensity functions are constant and the mean value formulas hold without reference to the points $z\in{\Bbb R}^d$. Beside the two mean value relations from above, general local mean value relation may be obtained by generalizing (\ref{EQMIXFK}) to the non-stationary case by a suitable concept of localization, see \cite[pp. 523-524]{SW}. However, the absence of scaling relations in the non-stationary case makes them less explicit.\\
We pass now to the stationary but not necessarily isotropic regime and assume for the remaining part of this subsection that the measure $\Lambda$ is translation-invariant rather than just translation-regular. In this case, the result of Theorem \ref{CORDISTR} remains also true, with $\PHT(s)$ replaced by $\PHT(s\Lambda)$ there. To see it, observe that the mean value formula $$N_{k,I}^{(d)}={d-k\over d}2^{d-k-1}\gamma_{\PHT(t\Lambda)}^{(k)}$$ is in force also in the stationary but non-isotropic case, in full analogy to (\ref{MIXTYPEQ}). Moreover, by translation invariance there exists some constant $c$, not depending of $t$, such that $\gamma_{\PHT(t\Lambda)}^{(k)}=ct^d$. Thus $${1\over s}{\gamma_{\PHT(s\Lambda)}^{(k)}\over N_{k,I}^{(d)}}={1\over s}{cs^d\over{d-k\over d}2^{d-k-1}ct^d}={d\over(d-k)2^{d-k-1}}{s^{d-1}\over t^d}$$ and we have again in analogy to (\ref{MIXTYPEQ})
\begin{eqnarray}
 \nonumber
{\Bbb Q}_k^{Y(t)} &=& (d-k) 2^{d-k-1} \int_0^t {1\over s}{\gamma_{\PHT(s\Lambda)}^{(k)}\over N_{k,I}^{(d)}}
 {\Bbb Q}_k^{\PHT(s\Lambda)} ds \\
 \nonumber &=& (d-k)2^{d-k-1}\int_0^t{d\over(d-k)2^{d-k-1}}{s^{d-1}\over t^d}{\Bbb Q}_k^{\PHT(s\Lambda)}ds
\end{eqnarray}
and hence 
\begin{equation}\label{formdirstat}
 {\Bbb Q}_k^{Y(t)} =\int_0^t{ds^{d-1}\over t^d}{\Bbb Q}_k^{\PHT(s\Lambda)}ds
\end{equation}
also in the non-isotropic case, where by ${\Bbb Q}_k^{\PHT(s\Lambda)}$ we understand the distribution of the typical $k$-dimensional I-face of the stationary Poisson hyperplane tessellation with intensity measure $s\Lambda$ and where $Y(t)=Y(t\Lambda)$ is the STIT tessellation with intensity measure $t\Lambda$. Also the time-marked analogue of (\ref{formdirstat}) holds true in the non-isotropic setting: \begin{equation}\hat{{\Bbb Q}}_k^{Y(t)}=\int_0^t{ds^{d-1}\over t^d}\left[{\Bbb Q}_k^{\PHT(s\Lambda)}\otimes\delta_s\right]ds.\label{formdirstattimemearked}\end{equation} Note that the statement about the mean f-vector of the typical $k$-dimensional I-face in Corollary \ref{corfvector} remains valid as well as the mean value formulas (\ref{NIKJ}), (\ref{MVF2}), (\ref{MVF23}), (\ref{MVF24}) and (\ref{MVF3}).\\ 
The equality (\ref{formdirstattimemearked}) from above allows us to determine the conditional distribution of the typical $k$-dimensional I-face of $Y(t)=Y(t\Lambda)$ given its birth time $0<s<t$. We readily see that this conditional distribution equals the distribution of the typical $k$-face of a Poisson hyperplane tessellation with the same parameter measure $t\Lambda$. Especially, the conditional length distribution of the typical I-segment given its birth time $s$ and direction $l$  is an exponential distribution with parameter $s\Lambda([e(l)])$ where $e(l)$ is a unit vector on $l,$ 
and thus the distribution of the length of the typical I-segment itself is a mixture of suitable exponential distributions. In particular, our considerations answer an extended version of the question formulated in the Outlook section of \cite{M09}. Moreover, the directional distribution ${\cal R}^{(k)}$ of the random process of $k$-dimensional I-faces, $1\leq k\leq d-2$, is in view of Corollary \ref{corDD} and (4.62,4.63) in \cite{SW} given by $${\cal R}^{(k)}(\cdot)={\kappa_{d-k}\over V_{d-k}(\Pi)}\rho_{(d-k)}^\perp(\cdot),$$ where by $\rho_{(d-k)}$ we denote the $(d-k)$-th projection generating measure of the \textit{associated zonoid} $\Pi$ \cite[(4.46)]{SW} -- an auxiliary convex body, which uniquely characterizes the law of the underlying STIT tessellation -- with generating hyperplane measure $t\Lambda,$
as defined by \cite[Eq. (14.36)]{SW}, whereas $(\cdot)^{\perp}$ stands for the orthogonal complement 
mapping.\\ Beside the distributional result for I-segments, mean values for the typical $k$-dimensional I-face of a stationary but non-isotropic STIT tessellation are also within reach of our methods and can be obtained by the same scheme as demonstrated above. In place of the concrete values involving Gamma-functions, in the non-isotropic case we have at our disposal the geometric parameters of the associated zonoid. As above, let $\Pi$ be the zonoid in the sense of \cite[Chap. 4.6]{SW} with translation-invariant generating hyperplane measure $t\Lambda$. Then the $j$-th intrinsic volume of the typical $k$-dimensional I-face is given by
\begin{eqnarray}
\nonumber {\Bbb E}V_j(I_k^{(d)}) &=& \int_0^t{ds^{d-1}\over t^d}{{d-j\choose d-k}V_{d-j}\left({s\over t}\Pi\right)\over{d\choose k}\Vol_d\left({s\over t}\Pi\right)}ds={{d-j\choose d-k}\over{d\choose k}}{V_{d-k}(\Pi)\over\Vol_d(\Pi)}\int_0^t{ds^{d-1}\over t^d}{\left({s\over t}\right)^{d-j}\over\left({s\over t}\right)^d}ds\\
 &=& {{d-j\choose d-k}\over{d\choose k}}{d\over d-j}{V_{d-k}(\Pi)\over\Vol_d(\Pi)}.\label{STATMEANVALUEFORMULA}
\end{eqnarray}
Here we have used (\ref{formdirstat}), \cite[Thm. 10.3.3]{SW} and the homogeneity of the intrinsic volumes. Note that (\ref{STATMEANVALUEFORMULA}) reduces to (\ref{eqEVJ}) in the isotropic case, i.e. when $\Lambda=\Lambda_{iso}$ or equivalently when the zonoid $\Pi$ is a $d$-dimensional ball with radius proportional to $t$.

\section{Second-Order Properties}\label{secSECOND}
 In this section we study second-order characteristics of iteration infinitely divisible random tessellations and stationary STIT tessellations. This is done from several different
 perspectives, both with the use of our martingale techniques and of recent very elegant second-order theory
 developed by Weiss, Ohser and Nagel \cite{NOW} for the planar case.

\subsection{Martingale Tools}\label{secSECONDMARTINGALETOOLS}
 The general martingale statements of Section \ref{secMARTINGALE} admit a convenient specialization to deal
 with the second-order characteristics of iteration infinitely divisible or stationary STIT tessellations as  well.
 Taking  $\phi$ to be a general bounded
 measurable functional of $(d-1)$-dimensional facets, regarded as usual as closed subsets of $W,$
 we  put $G(Y,t) := (\Sigma_{\phi}(Y) - {\Bbb E}\Sigma_{\phi}(Y(t,W)))^2$ so that
 $G(Y(t,W),t) = \bar\Sigma_{\phi}^2(Y(t,W))$ with 
 $\bar\Sigma_{\phi}(Y(t,W)) := \Sigma_{\phi}(Y(t,W)) - {\Bbb E}\Sigma_{\phi}(Y(t,W)),$
 and we use (\ref{EXPECT}) to check that
\begin{equation}\label{GPOCH}
 \frac{\partial}{\partial t} G(Y(t,W),t) = - 2 [\Sigma_{\phi}(Y(t,W)) - {\Bbb E}\Sigma_{\phi}(Y(t,W))]
 {\Bbb E}A_{\phi}(Y(t,W))
\end{equation}
where
\begin{equation}\label{APHI}
   A_{\phi}(Y) := \int_{[W]}\sum_{x\in\iota(H;Y)}\phi(\Facet(x;H|Y))\Lambda(dH).
\end{equation}
Put together (\ref{GEN}), (\ref{MART2}) and (\ref{GPOCH}) and use localization as in the discussion
preceding (\ref{EXPECT}) with $G_N,\; N \to \infty,$ chosen so that $(G_N(\cdot,\cdot) \wedge N) \vee -N 
\equiv (G(Y,t) \wedge N) \vee -N,$ that $|G_N(\cdot,\cdot)| \leq N+1$ and that $G_N(\cdot,t)$ be twice
continuously differentiable in $t,$ and with the localizing stopping times $T_N = \inf_{t \geq 0} 
(|G(Y(t,W),t)| \vee |\frac{\partial}{\partial t}G(Y(t,W),t)| \vee |\frac{\partial^2}{\partial t^2} G(Y(t,W),t)|)
 \geq N.$ Proceeding as there,  we readily conclude that
$$ \bar\Sigma_{\phi}^2(Y(t,W)) - \int_0^t \int_{[W]} \sum_{x\in\iota(H;Y(s,W))} \phi^2(\Facet(x;H|Y(s,W))) \Lambda(dH) ds + $$
$$ 2 \int_0^t [ \int_{[W]} \sum_{x\in\iota(H;Y(s,W))} \phi(\Facet(x;H|Y(s,W))) [\Sigma_{\phi}(Y(s,W))-{\Bbb E}\Sigma_{\phi}(Y(s,W))] \Lambda(dH) - $$
$$ [\Sigma_{\phi}(Y(s,W))-{\Bbb E}\Sigma_{\phi}(Y(s,W))] {\Bbb E} A_{\phi}(Y(s,W))] ds =  $$
\begin{equation}\label{VARMART}
 \bar\Sigma_{\phi}^2(Y(t,W)) - \int_0^t A_{\phi^2}(Y(s,W)) ds -
   2 \int_0^t \bar A_{\phi}(Y(s,W)) \bar\Sigma_{\phi}(Y(s,W)) ds 
\end{equation}
is a $\Im_t$-martingale as well, with $\bar A_{\phi}(Y(s,W)) := A_{\phi}(Y(s,W)) - {\Bbb E}A_{\phi}(Y(s,W)).$
Note for future reference that taking another bounded measurable facet functional $\psi,$ applying the above for
$\phi + \psi$ and $\phi -\psi$ and subtracting yields one further $\Im_t$-martingale
$$
 \bar\Sigma_{\phi}(Y(t,W)) \bar\Sigma_{\psi}(Y(t,W)) - \int_0^t A_{\phi \psi}(Y(s,W)) ds -  $$
\begin{equation}\label{VARMARTfutureref}
 \int_0^t (\bar A_{\phi}(Y(s,W)) \bar\Sigma_{\psi}(Y(s,W)) + \bar A_{\psi}(Y(s,W)) \bar\Sigma_{\phi}(Y(s,W))) ds.
\end{equation}
For general $\phi$ this cannot be simplified any further. However, in our further considerations
we shall focus our attention on  translation-invariant face functionals $\phi$ of the form
\begin{equation}\label{PHIFORM}
 \phi(f) := \Vol_{d-1}(f) \zeta(\vec{\bf n}(f))
\end{equation}
with $\vec{\bf n}(f)$ standing for the unit normal to $f$ and $\zeta$ for a bounded measurable
function on ${\cal S}_{d-1}.$ Then, using  (\ref{PHIFORM}) we see that 
$$ A_{\phi} \equiv \int_{[W]} \Vol_{d-1}(H \cap W) \zeta(\vec{\bf n}(H)) \Lambda(dH) = \const $$ 
and so $\bar A_{\phi} \equiv 0$ and thus, by (\ref{EXPECT}) and (\ref{VARMART}),
\begin{equation}\label{MART3}
  \bar\Sigma_{\phi}(Y(t,W)) \; \mbox{ and } \; \bar\Sigma_{\phi}^2(Y(t,W)) - \int_0^t A_{\phi^2}(Y(s,W)) ds
\end{equation}
are both $\Im_t$-martingales. In particular, see \cite[Thm. 4.2]{JS}, the martingale $\bar\Sigma_{\phi}(Y(t,W))$ has its
predictable quadratic variation process $\langle \bar\Sigma_{\phi}(Y(\cdot,W)) \rangle$ absolutely continuous and
given by
\begin{equation}\label{PREDQUVAR}
 \langle \bar\Sigma_{\phi}(Y(\cdot,W)) \rangle_t = \int_0^t A_{\phi^2}(Y(s,W)) ds.
\end{equation}
These observations are going to be crucial for our second-order analysis of
iteration infinitely divisible and stationary STIT tessellations, given in the next Subsection \ref{subsecVAR}.

\subsection{Variance of the Total Surface Area in a Window $W$}\label{subsecVAR}

It is the main purpose of this subsection to calculate the variance of the total surface area of the iteration infinitely divisible random tessellation $Y(t,W)$, for a compact convex window $W\subset{\Bbb R}^d$ with the property that $\Vol_d(W)>0$. We let $\Lambda$ be an arbitrary diffuse and locally finite measure on $\cal H$. Recall now (\ref{APHI}) and note that it implies
\begin{eqnarray}
\nonumber & & A_{\phi^2}(Y(t,W))=\int_{[W]} \sum_{f \in \Cells(Y \cap H)} \phi^2(f) \Lambda(dH)\\
&=&  \int_{[W]} \zeta^2(\vec{\bf n}(H)) \int_{W \cap H} \int_{W \cap H}
  {\bf 1}[x,y \mbox{ are in the same cell of } Y \cap H] dx dy 
  \Lambda(dH).\label{PHIKW}
\end{eqnarray} 
Thus, using (\ref{MART3}) and taking expectations of both sides yields immediately
 $$ \Var \Sigma_{\phi}(Y(t,W)) = $$
\begin{equation}\label{VAROG1}
  \int_0^t \int_{[W]} 
    \zeta^2(\vec{\bf n}(H)) \int_{H \cap W} \int_{H \cap W} 
    {\Bbb P}(x,y \mbox{ are in the same cell of } Y(s,W) \cap H) dx dy
    \Lambda(dH)ds. 
\end{equation}
Taking into account that $${\Bbb P}(x,y \mbox{ are in the same cell of } Y(s,W) \cap H) = \exp(-s\Lambda([xy])),$$ which follows from key property (a) and using (\ref{VAROG1}) we end up with
$$
 \Var(\Sigma_{\phi}(Y(t,W))) = \int_0^t \int_{[W]} 
    \zeta^2(\vec{\bf n}(H)) \int_{W \cap H} \int_{W \cap H} 
    \exp(-s\Lambda([xy])) dx dy \Lambda(dH) ds
$$
\begin{equation}\label{VAROG2}  
= \int_{[W]} \zeta^2(\vec{\bf n}(H)) \int_{W \cap H} \int_{W \cap H} 
    \frac{1-\exp(-t\Lambda([xy]))}{\Lambda([xy])} dx dy \Lambda(dH).\end{equation}For the stationary and isotropic case $\Lambda = \Lambda_{iso}$ we want to evaluate this integral further in the special case $\phi=\Vol_{d-1}$, i.e. when $\zeta\equiv 1$. First, we use the affine Blaschke-Petkantschin formula \cite[Thm. 7.2.7]{SW} with $q=1$ to obtain for any non-negative measurable function $h:({\Bbb R}^d)^2\rightarrow{\Bbb R}$ $$\int_{{\Bbb R}^d}\int_{{\Bbb R}^d}h(x,y)dxdy={d\kappa_d\over 2}\int_{{\cal L}}\int_L\int_L h(x,y)\left\|x-y\right\|^{d-1}\ell_L(dx)\ell_L(dy)dL,$$ where ${\cal L}$ is the space of lines in ${\Bbb R}^d$ with invariant measure $dL$, i.e. the affine $1$-dimensional Grassmannian in ${\Bbb R}^d$ and $\ell_L$ is the the Lebesgue measure on $L$ with normalization as specified in \cite[Thm. 13.2.12]{SW}. Taking now $$h(x,y)={\bf 1}[x\in W]{\bf 1}[y\in W]\left\|x-y\right\|^kg(x,y)$$ for some $k>-d$ and another non-negative measurable function $g:({\Bbb R}^d)^2\rightarrow{\Bbb R}^d$ we obtain
\begin{equation}
\int_W\int_W{\left\|x-y\right\|^kg(x,y)}dxdy={d\kappa_d\over 2}\int_{{\cal L}}\int_{W\cap L}\int_{W\cap L}\left\|x-y\right\|^{d-1+k}g(x,y)\ell_L(dx)d\ell_L(dy)dL.\label{EQBPF}
\end{equation}
For $k=-1$ this yields
\begin{equation}
\int_W\int_W{g(x,y)\over\left\|x-y\right\|}dxdy={d\kappa_d\over 2}\int_{{\cal L}}\int_{W\cap L}\int_{W\cap L}\left\|x-y\right\|^{d-2}g(x,y)\ell_L(dx)\ell_L(dy)dL.\label{CALC4}
\end{equation}
We replace now in (\ref{EQBPF}) for $k=0$, $W$ by $W\cap H$ for some fixed hyperplane $H$ and $d$ by $d-1$ and get
\begin{eqnarray}
\nonumber & & \int_{W\cap H}\int_{W\cap H}g(x,y)dxdy\\
\nonumber &=& {(d-1)\kappa_{d-1}\over 2}\int_{{\cal L}^H}\int_{W\cap H\cap L}\int_{W\cap H\cap L}\left\|x-y\right\|^{d-2}g(x,y)\ell_L(dx)\ell_L(dy)dL^H,
\end{eqnarray}
where by ${\cal L}^H$ we mean the $1$-dimensional affine Grassmannian restricted to $H$ with invariant measure $dL^H$. Averaging the last expression over all hyperplanes $H$ and using the fact that $\Lambda_{iso}(dH)\otimes dL^H=dL,$ see \cite[Thm. 13.2.12]{SW}, yields \begin{eqnarray}
\nonumber & & \int_{{\cal H}}\int_{W\cap H}\int_{W\cap H}g(x,y)dxdy\Lambda_{iso}(dH)\\
\nonumber &=& {(d-1)\kappa_{d-1}\over 2}\int_{{\cal H}}\int_{{\cal L}^H}\int_{W\cap H\cap L}\int_{W\cap H\cap L}\left\|x-y\right\|^{d-2}g(x,y)\ell_L(dx)\ell_L(dy)dL^H\Lambda_{iso}(dH)\\
&=& {(d-1)\kappa_{d-1}\over 2}\int_{{\cal L}}\int_{W\cap L}\int_{W\cap L}\left\|x-y\right\|^{d-2}g(x,y)\ell_L(dx)\ell_L(dy)dL.\label{CALC3}
\end{eqnarray}
By comparing (\ref{CALC4}) and (\ref{CALC3}) we finally conclude the non-trivial identity
\begin{equation}
\int_{[W]}\int_{W\cap H}\int_{W\cap H}g(x,y)dxdy\Lambda_{iso}(dH)={(d-1)\kappa_{d-1}\over d\kappa_d}\int_W\int_W{g(x,y)\over\left\|x-y\right\|}dxdy\label{CALC1}
\end{equation}
for any non-negative measurable function $g:({\Bbb R}^d)^2\rightarrow{\Bbb R}$. Using Equation (14) in \cite{NW05} together with the the mean projection formula \cite{SW}, Thm. 6.2.2 for $q=j=d-1$ (or
equivalently using the sectional property (b) in Subsection \ref{secSTIT}) we get from (\ref{CALC1}) with $$g(x,y)=\frac{1-\exp(-t\Lambda_{iso}([xy]))}{\Lambda_{iso}([xy])}={1-e^{-{2\kappa_{d-1}\over d\kappa_d}t\left\|x-y\right\|}\over {2\kappa_{d-1}\over d\kappa_d}\left\|x-y\right\|}$$ 
the following formula for $\Var(\Sigma_{\Vol_{d-1}}(Y(t,W)))=\Var(\Vol_{d-1}(Y(t,W)))$:
\begin{eqnarray}
\nonumber & & \int_{[W]}\int_{H\cap W}\int_{H\cap W}g(x,y)dxdy\Lambda_{iso}(dH)={d-1\over 2}\int_W\int_W{1-e^{-{2\kappa_{d-1}\over d\kappa_d}t\left\|x-y\right\|}\over\left\|x-y\right\|^2}dxdy\\
\nonumber &=& {d(d-1)\kappa_d\over 2}\int_0^\infty\overline{\gamma}_W(r){1-e^{-{2\kappa_{d-1}\over d\kappa_d}tr}\over r^2}r^{d-1}dr\\
\nonumber &=& {d(d-1)\kappa_d\over 2}\int_0^\infty\overline{\gamma}_W(r)r^{d-3}\left(1-e^{-{2\kappa_{d-1}\over d\kappa_d}tr}\right)dr
\end{eqnarray}
by using $d$-dimensional spherical coordinates. Here $$\overline{\gamma}_W(r)=\int_{{\cal S}_{d-1}}\Vol_d(W\cap(W+re_\varphi))\nu_{d-1}(d\varphi)$$ is the isotropized set-covariance function of the window $W$. Summarizing, we arrive at
\begin{theorem}\label{thmVAR} For the  stationary and isotropic STIT tessellation $Y(t)$ with surface intensity $t>0$ we have 
\begin{equation}
\Var(\Vol_{d-1}(Y(t,W)))= \frac{d-1}{2} 
 \int_W \int_W {1-e^{-{2\kappa_{d-1}\over d\kappa_d}t\left\|x-y\right\|}\over\left\|x-y\right\|^2}dxdy = 
\label{EQVAR0} \end{equation}
\begin{equation}
{d(d-1)\kappa_d\over 2}\int_0^\infty\overline{\gamma}_W(r)r^{d-3}\left(1-e^{-{2\kappa_{d-1}\over d\kappa_d}tr}\right)dr.\label{EQVAR}
\end{equation}
\end{theorem}
We specialize this now by taking $W$ to be the ball in ${\Bbb R}^d$ with radius $R>0$, i.e. $W=B_R^d$. In this case, the \textit{isotropized set-covariance function} $\overline{\gamma}_{B_R^d}(r)$ is given by $$\overline{\gamma}_{B_R^d}(r)=2R^d\kappa_{d-1}\int_{{r\over 2R}}^1(1-u^2)^{d-1\over 2}du=2\kappa_dR^d\left({1\over 2}-{r\over 2R}{\kappa_{d-1}\over\kappa_d}{_2F_1}\left({1\over 2},{1-d\over 2};{3\over 2};{r^2\over 4R^2}\right)\right)$$ (with $_2F_1$ being the Gauss hypergeometric function) for $0\leq r\leq 2R$ and $0$ otherwise, see \cite[Chap. 4.8.4]{HaSto}, and we have $$\Var(\Vol_{d-1}(Y(t,B_R^d)))={d(d-1)\kappa_d\kappa_{d-1}}R^d\int_0^{2R}\left(1-e^{-{2\kappa_{d-1}\over d\kappa_d}tr}\right)\int_{r\over 2R}^1(1-u^2)^{d-1\over 2}dur^{d-3}dr$$ $$={d(d-1)\kappa_d\kappa_{d-1}}R^d\int_0^{2R}\left(1-e^{-{2\kappa_{d-1}\over d\kappa_d}tr}\right)\left({1\over 2}-{r\over 2R}{\kappa_{d-1}\over\kappa_d}{_2F_1}\left({1\over 2},{1-d\over 2};{3\over 2};{r^2\over 4R^2}\right)\right)r^{d-3}dr,$$ where the constant before the integral may also be written as $${d(d-1)\kappa_d\kappa_{d-1}}R^d=(d-1){2\sqrt{\pi}\Gamma\left({d+1\over 2}\right)\over\Gamma\left({d\over 2}\right)}R^d={2(d-1)\over\lambda_1}R^d$$ with $\lambda_1$ from (\ref{LAMBDAK}). In the even more special case $d=3$ the isotropized set-covariance function $\overline{\gamma}_{B_R^3}(r)$ takes the form $$\overline{\gamma}_{B_R^3}(r)=\begin{cases}{4\pi\over 3}R^3\left(1-{3r\over 4R}+{r^3\over 16R^3}\right) &: 0\leq r\leq 2R\\ 0 &: r>2R\end{cases}$$
and the variance integral can be evaluated in a closed form: \begin{equation}\Var(\Vol_{2}(Y(t,B_R^3)))={4\pi^2\over 3t^4}\left(t^2R^2(12-8tR+3t^2R^2)+24(1+tR)e^{-tR}-24\right).\label{EXPLICITVAR3D}\end{equation} The same closed form cannot be obtained for $d=2$, since $\overline{\gamma}_{B_R^2}(r)$ has a more complicated structure, i.e. $$\overline{\gamma}_{B_R^2}(r)=2R^2\arccos\left({r\over 2R}\right)-{r\over 2}\sqrt{4R^2-r^2}$$ for $r$ between $0$ and $2R$ and $\overline{\gamma}_{B_R^2}(r)=0$ for $r>2R$. Unfortunately, the resulting integral can in this case not further be simplified.\\
Another important task in this context is to determine the large $R$ fixed $t$ asymptotics of the variance
 $\Var(\Vol_{d-1}(Y(t),W_R))$
for the family of growing windows $W_R=R\cdot W,\; R \to \infty.$ For $d=2$ we claim that 
\begin{equation}\label{EQVARASYM2D}
 \Var(\Vol_{1}(Y(t,W_R)) \sim \pi\Vol_2(W)R^2\log R,
\end{equation}
where $\sim$ stands for the asymptotic equivalence of functions, i.e. $f(R)\sim g(R)$ iff $f(R)/g(R)\rightarrow 1$ as $R\rightarrow\infty$. Indeed, this can be established by using (\ref{EQVAR}), the relation $\overline{\gamma}_{W_R} \sim\Vol_2(W_R)=R^2\Vol_2(W)$
valid uniformly for argument $r = O(R / \log R),$ the observation that $\overline{\gamma}_{W_R} \to 0$ for $r = \Omega(R \log R),$ 
together with the fact that $\int_0^{L(R)}(1-e^{-cr}){dr\over r}\sim \log R$, $c>0$, as soon as $\log L(R) \sim \log R,$
and the key property (c) of STIT tessellations: 
\begin{eqnarray}
\nonumber \Var(\Vol_1(Y(t,W_R))) &=& t^{-2}\Var(\Vol_1(Y(1,W_{tR}))) = {\pi\over t^2}\int_0^\infty\overline{\gamma}_{W_{tR}}(1-e^{-{2\over\pi}r}){dr\over r}\\
\nonumber &\sim& \pi t^{-2}\Vol_2(W_{tR})\log(tR)=\pi\Vol_2(W)R^2(\log R+\log t)\\
\nonumber &\sim& \pi R^2\Vol_2(W)\log R.
\end{eqnarray}
Another way to see it is to combine (\ref{VPHI}) below with the results from Subsection \ref{CLTplane}. Thus, as we see, in the planar case the studied asymptotics only depends on the area of $W.$ Things
get more complicated for $d > 2$ though. To see it, use (\ref{EQVAR0}) and the scaling property (c) of
STIT tessellations to obtain
\begin{eqnarray}
\nonumber \Var(\Vol_{d-1}(Y(t,W_R))) &=& R^{2(d-1)} \Var(\Vol_{d-1}(Y(Rt,W))\\
\nonumber &=& R^{2(d-1)}\frac{d-1}{2} 
    \int_W \int_W {1-e^{-{2\kappa_{d-1}\over d\kappa_d}Rt\left\|x-y\right\|}\over\left\|x-y\right\|^2}dxdy.
\end{eqnarray}
Consequently, we get for $d>2$ and $R \to \infty$
\begin{equation}\label{EQVARASYMHD}
\Var(\Vol_{d-1}(Y(t,W_R)))\sim R^{2(d-1)}{d-1\over 2}E_2(W)
\end{equation}
where $E_2(W)$ is the \textit{2-energy of $W$}, see \cite[Chap. 8]{MATI} given by
\begin{equation}\label{EN2}
 E_2(W) = \int_W \int_W \left\|x-y\right\|^{-2} dx dy . 
\end{equation}
Observe that this does not extend for the separately treated case $d=2$ because there the integral in (\ref{EN2}) diverges. It is easily seen that $E_2(\cdot)$ enjoys a \textit{superadditivity} property
\begin{equation}\label{SUPERADD}
 E_2(W_1 \cup W_2) \geq E_2(W_1) + E_2(W_2),\;\; W_1 \cap W_2 = \emptyset
\end{equation}
which stands in contrast to (\ref{EQVARASYM2D}) where the asymptotic expression is linear in $\Vol_2(W).$
We will now derive an integral geometric expression for this energy functional. Taking $g(x,y)\equiv 1$ and $k=-2$ in (\ref{EQBPF}) yields the remarkable identity
\begin{eqnarray}
\nonumber E_2(W) &=& \int_W\int_W{\left\|x-y\right\|^{-2}}dxdy={d\kappa_d\over 2}\int_{\cal L}\int_{W\cap L}\int_{W\cap L}\left\|x-y\right\|^{d-3}\ell_L(dx)\ell_L(dy)dL\\
 &=& {d\kappa_d\over(d-1)(d-2)}\int_{\cal L}\Vol_1(W\cap L)^{d-1}dL={2\over(d-1)(d-2)}I_{d-1}(W),\label{eee3comp}
\end{eqnarray}
with $I_{d-1}(W)$ being the $(d-1)$-st \textit{chord power integral} of $W$ in the sense of \cite[p. 363]{SW}. To display the dependency of the asymptotic variance $\Var(\Vol_{d-1}(Y(t,W_R)))$ on the geometry of $W$ we could take instead of the $(d-1)$-st chord power integral $I_{d-1}(W)$ of $W$ also the functional $$\int_{A(d,d-2)}V_{d-2}^2(W\cap E)dE,$$ where $A(d,d-2)$ is the set of affine $(d-2)$-planes in ${\Bbb R}^d$ with invariant measure $dE$, since $$I_{d-1}(W)={(d-1)d\kappa_d\over 2\kappa_{d-1}}\int_{A(d,d-2)}V_{d-2}^2(W\cap E)dE$$ according to \cite[Eq. (8.57)]{SW}. Thus, comparison with (\ref{eee3comp}) yields the alternative representation $$E_2(W)={d\kappa_d\over(d-2)\kappa_{d-1}}\int_{A(d,d-2)}V_{d-2}^2(W\cap E)dE.$$ The choice of one of these functionals is more or less a matter of taste, but we will take $I_{d-1}(W)$ here, because it will allow us in Section \ref{secCOMPARISON} a better comparison with other tessellations models. Hence, combining (\ref{eee3comp}) with (\ref{EQVARASYMHD}) from above, we arrive for $d\geq 3,\; R \to \infty$ at \begin{equation}\Var(\Vol_{d-1}(Y(t,W_R)))\sim {1\over d-2}I_{d-1}(W_R)={1\over d-2}R^{2(d-1)}I_{d-1}(W).\label{EQVARASYMHD2}\end{equation} In general, $I_{d-1}(W)$ cannot further be evaluated. But for $W=B_1^d$ we have by applying \cite{SW}, Theorem 8.6.6 (with a corrected constant) $$I_{d-1}(B_1^d)={d2^{d-2}}{\kappa_d\kappa_{2d-2}\over\kappa_{d-1}}$$ and, thus, the 2-energy of the $d$-dimensional unit ball $B_1^d$ equals $$E_2(B_1^d)={d2^{d-1}\over (d-1)(d-2)}{\kappa_d\kappa_{2d-2}\over\kappa_{d-1}}={2\pi^d\over(d-1)(d-2)}\Gamma\left({d\over 2}\right)^{-2}.$$ In the particular case $d=3$ we obtain the value $E_2(B_1^3)=4\pi^2$, which agrees with the explicit variance formula (\ref{EXPLICITVAR3D}). Another case, where the chord power integral can be evaluated is the practically relevant case of the cube $C_a^3$ in ${\Bbb R}^3$ with edge length $a>0$. Here we have $$I_2(C_1^3)={5\pi\over 3}-{13\over 9}-2\sqrt{2}\pi+2\pi\ln 2-{1\over 3}\ln 2+4\sqrt{2}\arctan\sqrt{2}-16F\approx 3.7557$$ with $$F=\int_{\sqrt{2}}^{\sqrt{3}}{\arctan\sqrt{x^2-2}\over x}dx,$$ for $a=1$ and in general $I_2(C_a^3)=a^4I_2(C_1^3)$ which can easily be obtained from \cite[Eq. (15)]{Itoh} and the fact that $I_{d-1}$ is homogeneous of degree $2d-2$.
\begin{remark}
Note that in the $2$-dimensional case the asymptotic variance of the total edge length in $W_R$ may also be written as $$\Var(\Vol_1(Y(t,W_R)))\sim I_1(W)R^2\log R,$$ yielding some kind of consistency with the higher dimensional cases represented by (\ref{EQVARASYMHD2}). But the argument using 2-energies given for $d\geq 3$ above cannot be applied to prove this formula.
\end{remark}
\begin{remark}
In our situation, \cite[Thm. 8.6.5]{SW} for $$f(r)=f(\Vert x-y\Vert)={1\over\Vert x-y\Vert}g(x,y)$$ can be applied, to deduce that $\Var(\Vol_{d-1}(Y(t,W)))$ is maximal among all compact convex bodies $W$ with positive volume exactly for $d$-dimensional balls, i.e. when $W=B_R^d$.
\end{remark}
\begin{remark}
In particular, Theorem \ref{thmVAR} establishes weak long range dependence (see \cite{HS}) present in stationary and isotropic STIT tessellations $Y(t)$, since $${\Var(\Vol_{d-1}(Y(t,W_R)))\over\Vol_d(W_R)}\sim\begin{cases} c_2\log R\rightarrow\infty &: d=2\\ c_d E_2(W)\Vol_d(W)^{-1}R^{d-2}\rightarrow\infty &: d\geq 3,\end{cases}$$ as $R\rightarrow\infty$ with $c_d$ being space dimension dependent constants. Heuristically this can be explained by the geometry of the I-facets of $Y(t)$. It was shown above that the $(d-1)$-volume of the typical I-facet has first but no second moment (see \cite{MNW07} in the planar case or \cite{T09} for a much more detailed discussion). Thus, the I-facets of $Y(t)$ can be extremely large, but there are no 'full' planes in the tessellations. This is reflected in the planar case by the log-term in the asymptotic expression for the variance of the total edge length, whereas in higher dimensions by the non-additivity of 2-energies determining the asymptotics. We will compare this with other tessellation models at the end of this paper.
\end{remark}
\begin{remark}
In the planar case, for $W=B_1^2$ the unit ball and $t=1$, (\ref{EQVARASYM2D}) was independently established by Lothar Heinrich (personal communication) using a quite different method which allows also the derivation of (\ref{EQVARASYM2D}) in its general form in the planar case $d=2$. The main idea is to dilate $Y(t,B_R^2)$ -- regarded as a random closed set in the plane -- by a small ball with radius $\varepsilon$ and to use the capacity functional of $Y(t,B_R^2)$ (see \cite{NW05}) together with a generalized Steiner formula and a limiting argument. Unfortunately the asymptotic formula for the variance is not sufficient to establish a central limit theorem as will be shown below in Section \ref{secCLT}.
\end{remark}

\subsection{Pair-Correlation Function and an Alternative Approach in the Planar Case}\label{secALTERNATIVE}

The \textit{second moment measure} $\mu_d^{(2)}(\cdot\times\cdot)$ of the random surface area measure of the stationary and isotropic STIT tessellation $Y(t)\subset{\Bbb R}^d$ can be defined by $$\mu_d^{(2)}(A_1\times A_2):={\Bbb E}[\Vol_{d-1}(Y(t,A_1))\Vol_{d-1}(Y(t,A_2))],$$ where $A_1$ and $A_2$ are bounded Borel sets in ${\Bbb R}^d$. By translation-invariance of the surface area measure of $Y(t)$ it is sufficient to regard the the \textit{reduced second moment measure} ${\cal K}_d$ of the random surface area measure of $Y(t)$, see \cite{SKM}, and since we are also in the isotropic case it is even enough to consider the \textit{reduced second moment function} $K_{d}(r)$ given by $$K_{d}(r):={\cal K}_d(B_r^d),\ \ r\geq 0.$$ Then quantity $t\cdot K_d(r)$ can be interpreted as the mean surface area of faces of $Y(t)$ within a ball with radius $r$ centered at a typical point of $Y(t)$ (when the tessellation is regarded under a suitable Palm distribution). In the case where $K_d(r)$ is differentiable in $r$ we can consider the \textit{pair-correlation function} $g_d(r)$ given by $$g_d(r)={1\over d\kappa_d r^{d-1}}{dK_d(r)\over dr},\ \ r\geq 0.$$ It describes the expected length density of $Y(t)$ at a given
distance $r$ from a typical point of $Y(t)$.\\ From the variance formula in Theorem \ref{thmVAR} the following can be easily deduced:
\begin{corollary}\label{corPCF} The pair-correlation function $g_d(r)$ of the random surface area measure of the stationary and isotropic random STIT tessellation $Y(t)$ is given by $$g_d(r)=1+{d-1\over 2t^2r^2}\left(1-e^{-{2\kappa_{d-1}\over d\kappa_d}tr}\right).$$
\end{corollary}
Especially for $d=2$, $g_d(r)$ becomes $$g_2(r)=1+{1\over 2t^2r^2}\left(1-e^{-{2\over\pi}rt}\right),$$ which was independently obtained by V. Weiss, J. Ohser and W. Nagel and will be presented in the forthcoming paper \cite{NOW}. We want to briefly sketch an alternative approach to obtain the variance of the total edge length of $Y(t,W),\; d =2,$ by using their recent results. The main result of \cite{NOW} says that the second moment measure of the random length measure induced by $Y(t)$ coincides with the second moment measure of the length measure induced by a Boolean segment process whose primary grain has length density $p_l(x)$ given by (\ref{EQPLANARPL}), which is the density of the length distribution of the typical I-segment of $Y(t)$ in the planar case. For a Boolean segment process with length distribution function $L$ and mean $\overline{l}$ of the typical segment we know that its reduced second moment function $K_2(r)$ is given by $$K_2(r)=\pi r^2+{1\over L_A\overline{l}}\left(\int_0^r x^2L(dx)+\int_r^\infty(2xr-r^2)L(dx)\right),$$ see \cite{S83}. Here $L_A$ is the mean segment length per unit area. In our special case, with $\overline{l}={\pi\over t}$ and $L_A=t$ and in view of the fact that length density (\ref{EQPLANARPL}) is a mixture of densities of exponential distributions, where the mixing density equals $2s\over t^2$, by applying repeated integration by parts we are led to
\begin{equation}K_2(r)=\pi r^2+{\pi\over t^2}\left(\gamma+\ln\left({2\over\pi}tr\right)+\Ei\left({2\over\pi}tr\right)\right),\label{KR}\end{equation} with $\gamma$ being the Euler--Mascheroni constant and $\Ei$ the exponential integral $$\Ei(x)=\int_1^\infty{e^{-xy}\over y}dy.$$ Using now the definition of the pair-correlation function $g_2(r)$ and the fact that the derivative of $\Ei(x)$ equals $-{e^{-x}/x}$ we get the following formula for the pair-correlation function $g_2(r)$ of the random length measure of $Y(t)$:
\begin{equation}
\nonumber g_2(r)={1\over 2\pi r}{dK(r)\over dr}=1+{1\over 2r^2t^2}\left(1-e^{-{2\over\pi}rt}\right) \label{GR}
\end{equation}
and this is the same as we obtained above in Corollary \ref{corPCF} for $d=2$. Using the pair-correlation function calculated in this way together with the formula $$\Var(\Vol_1(Y(t,W)))=2\pi t^2\int_0^\infty\overline{\gamma}_W(r)g_2(r)rdr-t^2\Vol_2^2(W)$$ from \cite[p. 233]{SKM} we can also obtain (\ref{EQVAR}) for $d=2$. It should be emphasized though that our original approach developed in Subsection
\ref{subsecVAR} above yields information also on higher dimensional cases. For example we have for the spatial case $d=3$
\begin{eqnarray}
\nonumber g_3(r) &=& 1+{1\over t^2r^2}\left(1-e^{-{1\over 2}tr}\right),\\
\nonumber K_3(r) &=& {4\pi\over 3}r^3+{4\pi\over 3t^3}\left(3tr-6+6e^{-{2\over 2}tr}\right).
\end{eqnarray}
The general expression for $K_d(r)$ is rather complicated and involves again special functions and is for this reason omitted.

\section{Central Limit Problem}\label{secCLT}
Having studied the first- and second-order properties of iteration infinitely divisible and stationary STIT tessellations, we now pass to the central limit problem in the stationary (thus STIT) regime. This problem will be considered in two
closely related settings, interestingly leading to results of very different qualitative natures. First, in Subsection \ref{secBROWNIANCONV} we shall focus our interest on {\it residual length/surface increment processes} arising respectively as cumulative length or surface area of I-facets born after a certain fixed time in the MNW-construction. In this set-up we shall establish a functional central limit theorem with the limit process identified as a suitable time-change of the standard Wiener process. Next, we shall pass to the full length/surface process, taking into account also the I-facets born at the very initial {\it big bang} stages of the MNW-construction, as descriptively termed in \cite{MNW}. It turns out that whereas in dimension $2$ the functional central limit theorem and Brownian convergence is preserved as shown in Subsection \ref{CLTplane}, this is no more the case for dimensions $3$ and higher, where non-Gaussian limits arise as argued in Subsection \ref{CLTspace}. This apparently surprising phenomenon is in fact due to the influence of the {\it big bang} phase itself, which is negligible in two dimensions but turns out crucial in higher dimensions. \\To proceed with our discussion, we put $W_R := R W$ for $R > 0$ and with  $W$ standing for a compact convex set of non-empty interior, to remain fixed throughout this section. We also let $\Lambda$ be some translation-invariant measure on $\cal H$ as in Section \ref{secSTIT}. Finally, we only shall consider $\phi$ of the form (\ref{PHIFORM}) in this section. 

\subsection{Brownian Convergence for Length/Surface Increment Processes}\label{secBROWNIANCONV}

 For fixed $s_0>0$ we consider the centered {\it surface increment process}
 $$ \left(\frac{1}{R^{d/2}} [\bar \Sigma_{\phi}(Y(t,W_R)) - \bar \Sigma_{\phi}(Y(s_0,W_R))] \right)_{t\in[s_0,1]} $$
 corresponding to the increments of $\Sigma_{\phi}(Y(t, W_R))$ relative to some initial time moment $s_0$ in the
 MNW-construction.  The main theorem of this subsection is
\begin{theorem}\label{BrownianConvergence}
 For each $s_0 > 0$ the centered surface increment process
 $$ \left( {\cal S}^{R,W}_{s_0,t} := 
     \frac{1}{R^{d/2}} [\bar \Sigma_{\phi}(Y(t,W_R)) - \bar \Sigma_{\phi}(Y(s_0,W_R))] \right)_{t\in[s_0,1]}, $$
 converges in law, as $R\to\infty,$ on the space ${\cal D}[s_0,1]$ of right continuous functions with left-hand limits
 (c\`adl\`ag)  on $[s_0,1]$ endowed with  
 the usual Skorokhod topology \cite[Chap. 3, Sec. 14]{BI}, to time-changed Wiener process
 $$ t \mapsto {\cal W}_{V_W(\phi,\Lambda) \int_{s_0}^t s^{1-d} ds}, $$
 where ${\cal W}_{(\cdot)}$ is the standard Wiener process and $V_W(\phi,\Lambda)$ is given
 by (\ref{VPHI}) or alternatively (\ref{VPHIALTERNATIVE}) below. In particular,
 $$ {\cal S}^{R,W}_{s_0,1} = \frac{1}{R^{d/2}} [\bar \Sigma_{\phi}(Y(1,W_R)) - \bar\Sigma_{\phi}(Y(s_0,W_R))] $$
 converges in law to ${\cal N}(0,V_W(\phi,\Lambda)\int_{s_0}^1 s^{1-d} ds)$, a normal distribution with mean $0$ and variance $V_W(\phi,\Lambda)\int_{s_0}^1 s^{1-d} ds$.
\end{theorem}
Note that this statement cannot be extended to  $s_0 \downarrow 0,$ as would be
of interest as potentially yielding Gaussian limit for $\bar \Sigma_{\phi}(Y(1,W_R)).$
The problem is that the variance integral $V_W(\phi,\Lambda) \int_{s_0}^t s^{1-d} ds$ diverges at
$0.$ As already signalled above, we will see below that this difficulty can be overcome for $d=2$ but not for $d>2$ where
we will obtain non-Gaussian limits for suitably normalized $\bar \Sigma_{\phi}(Y(1,W_R)).$

\paragraph{Proof of Theorem \ref{BrownianConvergence}}
Note first that 
$$
 \frac{1}{R^d} A_{\phi^2}(Y(1,W_R)) = \frac{1}{R} \int_{[W_R]} \frac{1}{R^{d-1}}
    \zeta^2(\vec{\bf n}(H)) \sum_{f \in \Cells(H \cap Y(1,W_R))} \Vol^2_{d-1}(f) \Lambda(dH)  = $$
\begin{equation}\label{NormVar}
 \int_{[W]} \frac{1}{R^{d-1}}
    \zeta^2(\vec{\bf n}(H)) \sum_{f \in \Cells(R H \cap Y(1,W_R))} \Vol^2_{d-1}(f) \Lambda(dH) .
\end{equation}
We claim that upon letting $R\to\infty$ this converges in probability to
$$
 V_W(\phi,\Lambda) :=
 \int_{[W]} \zeta^2(\vec{\bf n}(H)) \Vol_{d-1}(H \cap W) \frac{{\Bbb E}\Vol^2_{d-1}(\TypicalCell(H \cap Y(1)))}
   {{\Bbb E}\Vol_{d-1}(\TypicalCell(H \cap Y(1)))} \Lambda(dH)
$$
\begin{equation}\label{VPHI}
  = \Vol_d(W) \int_{{\cal S}_{d-1}} \zeta^2(u)  \frac{{\Bbb E}\Vol^2_{d-1}(\TypicalCell(u^{\bot} \cap Y(1)))}
   {{\Bbb E}\Vol_{d-1}(\TypicalCell(u^{\bot} \cap Y(1)))} {\cal R}(du),
\end{equation}
where $u^{\bot}$ is the orthogonal complement of $u\in{\cal S}_{d-1}$ and $\cal R$ is the directional distribution of the stationary STIT tessellations $Y(t)$ as given in (\ref{LADEF}). To see it, recall that $R H \cap Y(1)$ is a STIT
tessellation in $RH$ for each $R >0$ and $H \in {\cal H}.$ Thus, applying \cite[(4.6) and Thm 4.1.3]{SW} and (10.4) ibidem to this tessellation, we have $$\lim_{R\to\infty} \frac{1}{R^{d-1}} {\Bbb E} \sum_{f \in \Cells(R H \cap Y(1,W_R))} \Vol^2_{d-1}(f)$$
\begin{equation}\label{EXPE}
    =\Vol_{d-1}(H \cap W) \frac{{\Bbb E}\Vol^2_{d-1}(\TypicalCell(H \cap Y(1)))}
    {{\Bbb E}\Vol_{d-1}(\TypicalCell(H \cap Y(1)))}. 
\end{equation} 
Next, we observe that $R H \cap Y(1,W_R) \overset{D}{=} R \cdot_H (H \cap W) \cap Y(1)$ where
$\cdot_H$ is the scalar multiplication relative in $H,$ that is to say $H \ni R \cdot_H x = p_H(0) + R (x-p_H(0)),\; x \in H$
with $p_H$ standing for the orthogonal projection on $H.$  Thus, using the strong mixing and tail triviality theory for
STIT tessellations recently developed by Lachi\`eze-Rey \cite[Thm 2]{LR}, noting that tail trivial stationary processes are ergodic \cite[Prop. 14.9]{GEO} and then using the standard multidimensional ergodic theorem, see e.g. Cor. 14.A5 ibidem,
to $\frac{1}{R^{d-1}} \sum_{f \in \Cells(R \cdot_H (H \cap W) \cap Y(1))} \Vol^2_{d-1}(f),$ we get from (\ref{EXPE}) that
$$\lim_{R\to\infty} \frac{1}{R^{d-1}} \sum_{f \in \Cells(R H \cap Y(1,W_R))} \Vol^2_{d-1}(f)$$
\begin{equation}\label{EXPEL1}
    =\Vol_{d-1}(H \cap W) \frac{{\Bbb E}\Vol^2_{d-1}(\TypicalCell(H \cap Y(1)))}
    {{\Bbb E}\Vol_{d-1}(\TypicalCell(H \cap Y(1)))} 
\end{equation} 
in probability. Putting this together with (\ref{NormVar}) and integrating over $[W]$ yields
\begin{equation}\label{PRZEDZBIEZN}
 \lim_{R \to \infty} \frac{1}{R^d} A_{\phi^2}(Y(1,W_R)) = V_W(\phi,\Lambda) \;\; \mbox{ in probability } 
\end{equation}
as required. Note now that by scaling properties of $Y(s,W_R)$ (see key property (c)) and $\phi^2$ for $s > 0$ we have
\begin{equation}\label{SCALINGIDS}
 \frac{1}{R^d} A_{\phi^2}(Y(s,W_R)) \overset{D}{=} \frac{1}{R^d s^{2d-1}} A_{\phi^2}(Y(1,W_{sR})) \overset{D}{=}
 \frac{1}{s^{d-1}} \frac{1}{(Rs)^d} A_{\phi^2}(Y(1,W_{sR})).
\end{equation}
Thus, combining (\ref{PRZEDZBIEZN}) with the scaling relation (\ref{SCALINGIDS}) we get
\begin{equation}\label{CONVERGSTAT}
 \lim_{R\to\infty} \frac{1}{R^d} A_{\phi^2}(Y(s,W_R)) = \frac{1}{s^{d-1}} V_W(\phi,\Lambda)\;\; \mbox{ in probability uniformly in } s \in [s_0,1].
\end{equation}
This crucial statement puts us now in context of the general martingale limit theory. Indeed, using (\ref{MART3}) we
see that ${\cal S}^{R,W}_{s_0,s} = \frac{1}{R^{d/2}} [\bar \Sigma_{\phi}(Y(1,W_R)) - \bar\Sigma_{\phi}(Y(s_0,W_R))]$ is a martingale with absolutely continuous predictable quadratic variation process
\begin{equation}\label{QUADVAR}
 \langle {\cal S}^{R,W}_{s_0,\cdot} \rangle_t = \int_{s_0}^t \frac{1}{R^d} A_{\phi^2}(Y(s,W_R)) ds,
\end{equation}
 see \cite[Thm. 4.2]{JS}.
In these terms, (\ref{CONVERGSTAT}) yields for each $t$ 
\begin{equation}\label{QVARCONVREL}
 \lim_{R \to \infty} \langle {\cal S}^{R,W}_{s_0,\cdot} \rangle_t  = \int_{s_0}^t \frac{1}{s^{d-1}} V_W(\phi,\Lambda)
 \;\; \mbox{ in probability. }
\end{equation}
 We want now to apply the martingale functional limit theorem, see e.g. \cite[Thm 2.1]{WW}. The condition
 (ii.6) there is just (\ref{QVARCONVREL}) whereas condition (ii.4) there is trivially verified because the predictable quadratic variation $\langle {\cal S}^{R,W}_{s_0,\cdot} \rangle$ has no jumps by (\ref{QUADVAR}).
It remains to check the condition (ii.5) ibidem, which is that the second moment of the maximum jump
${\cal J}({\cal S}^{R,W}_{s_0,\cdot};1)$ of the process $({\cal S}^{R,W}_{s_0,s})_{s \in [s_0,1]}$
 goes to $0$ as $R \to \infty.$ To this end, note first that, with probability one,
 ${\cal J}({\cal S}^{R,W}_{s_0,\cdot};1)$  is bounded above by a constant multiple of
 $R^{-d/2}$ times the $(d-1)$-th
 power of the diameter of the largest cell of $Y(s_0,W_R).$ Since the typical cell of $Y(s_0)$
 is Poisson with intensity measure $s_0 \Lambda$ by property (a) of STIT tessellations,
 we conclude by standard properties that the expected number of cells in $Y(s_0,W_R)$ 
 with diameters exceeding $D$ is of the order $O(R^d \exp(-D))$ since  $\Lambda$ has been
 assumed to have its support spanning the whole of ${\Bbb R}^d,$ see Subsection \ref{secSTIT}.
 Summarizing, setting $u = D^{d-1} R^{-d/2}$ we are led to
 \begin{equation}\label{JBOUND}
  {\Bbb P}({\cal J}({\cal S}^{R,W}_{s_0,\cdot};1) > u) = O(R^d \exp(-R^{d/(2d-2)} u^{1/(d-1)})).
 \end{equation}
 Clearly, (\ref{JBOUND}) is much more than enough to guarantee that  
 $$ \lim_{R\to\infty}  {\Bbb E}{\cal J}^2({\cal S}^{R,W}_{s_0,\cdot};1) = 0$$
 which gives the required condition (ii.5) of Theorem 2.1 in \cite{WW}. 
Upon a trivial time change, this theorem yields now the functional convergence
in law as stated in our Theorem \ref{BrownianConvergence}. This completes the proof. $\hfill\Box$ \\ \\ An alternative formula can be provided for the factor $V_W(\phi,\Lambda)$. Denote by $\Pi$ the associated zonoid of STIT tessellation with generating hyperplane measure $\Lambda$, by $\Pi^o$ its \textit{dual body} and by ${\cal R}$ the directional distribution of the STIT tessellation. Then we have 
\begin{proposition}\label{propVW} It holds \begin{equation}V_W(\phi,\Lambda)=\Vol_{d}(W){(d-1)!\over 2^{d-1}}\int_{{\cal S}_{d-1}}\zeta^2(u)\Vol_{d-1}((\Pi|u^\perp)^o){\cal R}(du),\label{VPHIALTERNATIVE}\end{equation}
where $\Pi|u^\perp$ stands for the orthogonal projection of $\Pi$ onto the hyperplane $u^\perp$ and where the
polar body $(\Pi|u^{\perp})^o$ is considered relative to $u^{\perp}.$ In the isotropic case, i.e. when ${\cal R}=\nu_{d-1}$, this reduces to $$V_W(\phi,\Lambda_{iso})=\Vol_d(W)2^{d-1}\pi^{d-{3\over 2}}\Gamma\left({d+1\over 2}\right)^{d-1}\Gamma\left({d\over 2}\right)^{2-d}\int_{{\cal S}_{d-1}}\zeta^2(u)\nu_{d-1}(du).$$
\end{proposition}
In particular for $\zeta\equiv 1$, $W=B_1^d$ the unit ball and $d=2$ and $d=3$ we get the values $$V_{B_1^2}(\Vol_1,\Lambda_{iso})=\pi^2\ \ \ \ \ \text{and}\ \ \ \ V_{B_1^3}(\Vol_2,\Lambda_{iso})={32\over 3}\pi^2,$$ respectively.
\paragraph{Proof of Proposition \ref{propVW}} First, \cite[Cor. 3.7]{FW} provides a formula for the second moment of the volume of the typical Poisson cell of a stationary Poisson hyperplane tessellation in ${\Bbb R}^d$. In terms of the zonoid $\Pi$ it reads $${\Bbb E}\Vol_d^2(\TypicalCell(\PHT(\Lambda)))={d!\over 2^d}{\Vol_d(\Pi^o)\over\Vol_d(\Pi)},$$ where we have used formula \cite[(4.63)]{SW}. Moreover, the mean volume of $\TypicalCell(\PHT(\Lambda))$ is given by $${\Bbb E}\Vol_d(\TypicalCell(\PHT(\Lambda)))={1\over\Vol_d(\Pi)}$$ according to \cite[Thm. 10.3.3 and (10.4)]{SW}. Using now Eq. (4.61) ibidem, key property (a) of STIT tessellations and replacing $d$ by $d-1$ in the last to formulas we obtain (\ref{VPHIALTERNATIVE}).\hfill $\Box$

\subsection{Gaussian Limits in the Planar Case}\label{CLTplane}
In the planar case the asymptotic behaviour of $\bar \Sigma_{\phi}(Y(1,W_R))$ turns out to be Gaussian.
We write
\begin{equation}\label{TAUDEF}
 \tau(s,R) := \exp([\log R - \log\log R](s-1)) = R^{s-1} (\log R)^{1-s}
\end{equation}
and we define the total length process
\begin{equation}\label{LGTHPROC}
 {\cal L}^{R,W}_s := \frac{1}{R \sqrt{\log R}} \bar\Sigma_{\phi}(Y(\tau(s,R),W_R)),\; s \in [0,1].
\end{equation}
 The main result of this subsection is 
 \begin{theorem}\label{CLT2D}
  The total length process $({\cal L}^{R,W}_s)_{s\in [0,1]}$ converges in law, as $R\to\infty,$ on the space
  ${\cal D}[0,1]$ of c\`adl\`ag functions on $[0,1]$ endowed with the usual Skorokhod topology, to 
 $(\sqrt{V_W(\phi,\Lambda)} {\cal W}_s)_{s\in [0,1]}$ where, again, ${\cal W}_{(\cdot)}$ stands for the
 standard Wiener process.
\end{theorem}
In particular, for $\phi=\Vol_1$ we have 
$$\frac{1}{R \sqrt{\log R}} \bar\Sigma_{\Vol_1}(Y(1,W_R))\Longrightarrow{\cal N}(0,\pi\Vol_2(W))$$
for the stationary and isotropic STIT tessellation $Y(1)$ in the plane, where $\Longrightarrow$ means
convergence in law.

\paragraph{Proof of Theorem \ref{CLT2D}}
 Note first that
\begin{equation}\label{TAUPROP}
 \tau(0,R) = \frac{\log R}{R},\; \tau(1,R) = 1,\; \frac{\partial}{\partial s}\tau(s,R) = \tau(s,R) [\log R - \log\log R].
\end{equation}
 Thus, defining the auxiliary process
 \begin{equation}\label{PROCDEF}
  M^{R,W}_s = M_s := \frac{1}{R \sqrt{\log R - \log \log R}} [\bar\Sigma_{\phi}(Y(\tau(s,R),W_R)) - \bar\Sigma_{\phi}(Y(\tau(0,R),W_R))]
\end{equation}
 and using (\ref{MART3}) with $W_R := R W$ and under variable substitution $s := \tau(u,R)$ and $t := s$ 
 with LHS variables corresponding to the notation of (\ref{MART3}) and RHS to that used here,
 we see that, by (\ref{TAUPROP}),
\begin{equation}\label{MART4}
 (M_s)_{s=0}^1 \mbox{ and } \left( M^2_s - \int_0^s  \frac{\tau(u,R)}{R^2} A_{\phi^2}(Y(\tau(u,R), W_R)) du \right)_{s\in [0,1]}
\end{equation}
are $\Im_{\tau(s,R)}$-martingales. In particular, see \cite[Thm. 4.2]{JS}, the predictable quadratic variation process
$\langle M \rangle_s$ is given by
\begin{equation}\label{MQV}
 \langle M \rangle_s = \int_0^s \frac{\tau(u,R)}{R^2} A_{\phi^2}(Y(\tau(u,R),W_R)) du,\; s \in [0,1].
\end{equation} 
Repeating the argument leading to (\ref{CONVERGSTAT}) we see that
\begin{equation}\label{CONVSTAT2}
 \lim_{R\to\infty} \frac{\tau(s,R)}{R^2} A_{\phi^2}(Y(\tau(s,R), W_R)) = V_W(\phi,\Lambda) \;\;
    \mbox{ in probability, uniformly in } s \in [0,1]. 
\end{equation}
Note that the uniformity in $s$ comes, as in the case of (\ref{CONVERGSTAT}), from the relation (\ref{SCALINGIDS})
implying that, in distribution, all instances of the LHS for different values of $s$ are just scaling instances of the same
object $\tilde{R}^{-2} A_{\phi^2}(Y(1,W_{\tilde{R}}))$ for $\tilde{R} = R / \tau(s,R)$ and thus, in terms of the considered
convergence in probability to a deterministic limit, we are just dealing with a single asymptotic statement. Consequently, by (\ref{CONVSTAT2}) and in full analogy to (\ref{QVARCONVREL}), 
\begin{equation}\label{QVARCONVREL2}
 \lim_{R\to\infty}   \langle M \rangle_s =  \int_0^s  V_W(\phi,\Lambda) du = s V_W(\phi,\Lambda) \;\; \mbox{ in probability.}
\end{equation} 
Thus, we are again in a position to apply \cite[Thm 2.1]{WW} yielding the functional convergence in law, as $R \to \infty,$ in ${\cal D}[0,1]$ of $(M_s)_{s\in [0,1]}$ to $(\sqrt{V_W(\phi,\Lambda)} {\cal W}_s)_{s\in [0,1]}.$ Indeed,
condition (ii.6) there is just (\ref{QVARCONVREL2}), condition (ii.4) is trivial in view of (\ref{MQV}), whereas
the condition (ii.5) is verified by noting that, with probability one, ${\cal J}(M;1) = \frac{1}{R \sqrt{\log R}}
O(R \diam(W)) = O(1/\sqrt{\log R})$ so that in particular $\lim_{R\to\infty} {\Bbb E}{\cal J}^2(M;1) = 0$
as required.  Denoting now by
$C^{R,W}$ the {\it correction term} $(R \sqrt{\log R})^{-1} \bar\Sigma_{\phi}(Y(\tau(0,R),W_R))$ such that 
$$ {\cal L}^{R,W}_s = C^{R,W} + \sqrt{\frac{\log R - \log\log R}{\log R}} M_s, $$
noting that $\log R - \log\log R \sim \log R$ and that, by the scaling property (c) of STIT tessellations and
by (\ref{EQVARASYM2D}), 
\begin{equation}\label{NEGcorr}
 \Var(C^{W,R}) = O([R^{-2} (\log R)^{-1}] [R^2/(\log R)^2] [(\log R)^2 (\log\log R)])  = O(\log\log R/ \log R)
\end{equation}
we see that the processes $M_s$ and ${\cal L}^{R,W}_s$ are asymptotically equivalent in ${\cal D}[0,1]$ as $R\to\infty.$ 
This completes the proof of Theorem \ref{CLT2D}. $\hfill\Box$ 

\begin{remark}\label{CONVrem}
 In the context of proof of Theorem \ref{CLT2D} it should be remarked that the 'negligible correction term'
 $C^{R,W}$ has its variance of order $O(\log\log R / \log R)$ and thus indeed tending to $0,$
 but extremely slowly. Consequently, although the Gaussian CLT holds for ${\cal L}^{R,W}_0,$ it is quite natural
 to expect that the convergence rates are extremely slow, conjecturedly logarithmic. This is due to
 the fact that dimension $2$ is the largest dimension (critical dimension) where the Gaussian limits are still
 present. In dimensions $3$ and higher there is no Gaussian CLT and the 'correction term' analogous to
 $C^{R,W}$ will turn out order-determining rather than negligible, as shown in the next subsection.
\end{remark}

\subsection{Non-Gaussian Limits for $d > 2$}\label{CLTspace}
 We claim that the argument of Subsection \ref{CLTplane} above cannot be repeated for $d > 2.$
 Intuitively, this is due to the fact that for $d > 2$
the variance order of $\bar\Sigma_{\phi}(Y(1,W_R))$ is $O(R^{2(d-1)})$ and so even the very first faces
born in the cell division process already do bring a non-negligible contribution to the overall variance. Thus, we
cannot split the whole STIT construction into the {\it warm-up phase} ($t \in [0,R^{-1} \log R]$ for $d=2$)
with negligible variance contribution and the {\it proper phase} unfolding already in a typical STIT environment.
In fact, we claim CLT does not hold for STIT length functionals in dimension greater than $2!$ To see it,  observe
first that, by the scaling property (c) of STIT tessellations,
\begin{equation}\label{RNIEDOCTG}
 R^{-(d-1)} \bar\Sigma_{\phi}(Y(1,W_R))  \overset{D}{=} \bar\Sigma_{\phi}(Y(R,W)).
\end{equation}
Further, recall that by (\ref{MART3}) the process $R \mapsto \bar\Sigma_{\phi}(Y(R,W))$ is a square-integrable
martingale with absolutely continuous predictable quadratic variation process given in (\ref{PREDQUVAR})
and, moreover, by (\ref{VAROG2}) we have
$$ {\Bbb E}\bar\Sigma^2_{\phi}(Y(R,W)) = \int_{[W]} \zeta^2(\vec{\bf n}(H)) \int_{W \cap H}
    \int_{W \cap H} \frac{1-\exp(-R\Lambda([xy]))}{\Lambda([xy])} dx dy \Lambda(dH) $$
which is bounded uniformly in $R.$
Consequently, by the martingale convergence theorem, there exists a centered square-integrable
random variable $\Xi(W)$ such that
\begin{equation}\label{OSTCONV}
 \Xi(W) = \lim_{R\to\infty} \bar\Sigma_{\phi}(Y(R,W))
\end{equation}
 a.s. and in $L^2$ and, moreover, 
\begin{equation}\label{XIVAR}
 \Var\Xi(W) = \int_{[W]} \zeta^2(\vec{\bf n}(H)) \int_{W \cap H}
    \int_{W \cap H} \frac{1}{\Lambda([xy])} dx dy \Lambda(dH). 
\end{equation}
Using now (\ref{RNIEDOCTG}) we readily conclude that
\begin{equation}\label{CTGDW2}
 R^{-(d-1)} \bar\Sigma_{\phi}(Y(1,W)) \Longrightarrow \Xi(W)
\end{equation}
as $R \to \infty.$\\ We claim that the variable $\Xi(W)$ is not normal. Even though we are able to show this fact for all
$W$ and translation invariant $\Lambda$ by establishing non-Gaussian tail decay, for simplicity we
only give a proof for an easily tractable
particular case, postponing the study of more involved properties of the random field $\Xi(W),\;
W \subseteq {\Bbb R}^d$ for a future paper. Namely, take $W = [0,1]^d$ and 
\begin{equation}\label{PARTICULARLA}
  \Lambda := \sum_{i=1}^d \int_{-\infty}^{+\infty} \delta_{r e_i + e_i^{\bot}} dr
\end{equation}
where $e_i,\; i =1,\ldots,d$ are vectors of the standard orthonormal basis for ${\Bbb R}^d$ and
$\delta_{r e_i + e_i^{\bot}}$ is the unit mass concentrated on the hyperplane orthogonal to $e_i$
in distance $r$ from the origin. Consider the event ${\cal E}_{N},\; N >0,$ that only hyperplanes orthogonal
to $e_1$ have been born during the time $[0,1]$ of the MNW-construction and their number exceeds
$N.$ Observe that, in view of the form (\ref{PARTICULARLA}) of $\Lambda,$ 
$ {\Bbb P}({\cal E}_N) = \exp(-d) \sum_{k=N+1}^{\infty} \frac{1}{k!}  $
and thus
\begin{equation}\label{TAILSN}
 \log {\Bbb P}({\cal E}_N) = - \Theta(N \log N),
\end{equation} 
where by $\Theta(\cdot)$ we mean something bounded both from below and above by multiplicities of the argument, i.e. $\Theta(\cdot) = O(\cdot) \cap \Omega(\cdot).$ Further, given a fixed collection of all hyperplanes $H_1,\ldots,H_k,\; k > N,$ born at times
between $0$ and $1,$ on the event ${\cal E}_N$ we see that the conditional law of $\Xi(W)$ 
coincides with that of $k-d$ plus  sum of independent copies $\xi_1,\ldots,\xi_{k+1}$ of
$\Xi(W_1),\ldots,\Xi(W_{k+1})$ respectively, 
where $W_j,\; j=1,\ldots,k+1,$ are parallelepipeds into which $W$ is partitioned by $H_1,\ldots,H_k.$
Note that the extra $k$ above is the sum of $(d-1)$-volumes of $H_i \cap W$ whereas $-d = -
{\Bbb E}\Sigma_{\phi}(Y(1,W))$ is the centering term. Since $\Var[\xi_1+\ldots+\xi_{k+1}] =
\sum_{j=1}^{k+1} \Var\Xi(W_j)$ which is bounded above by $\Var\Xi(W)$ in view of (\ref{XIVAR}),
by Chebyshev's inequality we get ${\Bbb P}(\xi_1+\ldots+\xi_{k+1} \geq - 2 \sqrt{\Var\Xi(W)}) \geq 3/4.$
Thus, in view of (\ref{TAILSN})
\begin{equation}\label{TAILSXI}
 {\Bbb P}(\Xi(W) > N)  \geq \frac{3}{4} {\Bbb P}\left({\cal E}_{N + 2\sqrt{\Var\Xi(W)} + d}\right)
 = \exp(-\Theta(N\log N)).
\end{equation}
Since Gaussian variables exhibit tail decay of the order $\exp(-\Theta(N^2)),$ the random variable
$\Xi(W)$ cannot be normal, which completes our argument.

\section{Comparison with other Tessellation Models}\label{secCOMPARISON}

This section is devoted to a comparison of the stationary and isotropic STIT tessellations $Y(t)$ with stationary and isotropic Poisson hyperplane tessellations and Poisson-Voronoi tessellation with the same surface intensity $t>0$ stemming from a stationary Poisson process (we restrict ourselves to the isotropic case, since Poisson-Voronoi tessellations based on stationary Poisson processes are automatically isotropic). This is of particular interest for the statistical analysis of random tessellations, for example when constructing asymptotic confidence intervals for mean value estimators and for goodness-of-fit tests. A first example was considered in in \cite{NW082} where the linear contact distribution was used as a criterion. The mean values obtained in \cite{NW06}, \cite{NW08} and \cite{TW} can be used for such a test and now also some second order quantities. More details on goodness-of-fit tests for random tessellations can for example be found in \cite{HS} together with a couple of examples.\\ Denote by $\PHT(t,W)$ the stationary and isotropic Poisson hyperplane tessellation inside a bounded convex window $W\subset{\Bbb R}^d$. For $d=2$ we will write $\PLT(t,W)$ instead of $\PHT(t,W)$, since hyperplanes are just lines in this case. It was shown in \cite{H} that the asymptotic variance of the total $k$-volume of the $k$-faces of $\PHT(t,W)$ equals $$\Var(\Vol_k(\PHT(t,W)))\sim t^{2(d-k)-1}(d-k)^2{d\choose k}^2{\kappa_d^2\over\kappa_k^2}\left({\kappa_{d-1}\over d\kappa_d}\right)^{2(d-k)}J(W)$$ with \begin{equation}J(W)=\int_{{\cal H}}\Vol_{d-1}^2(W\cap H)\Lambda_{iso}(dH).\label{EQJD}\end{equation} In the case $W=B_R^d$ this integral can be evaluated explicitly and we have $$J(B_R^d)={((d-1)!\kappa_{d-1})^2 \over(2d-1)!}(2R)^{2d-1}.$$ For Poisson-Voronoi tessellations $\PVT(\gamma,W)$ (here restricted to $W$) constructed from a stationary Poisson point processes in ${\Bbb R}^d$ with intensity $\gamma$ we have from \cite{HSS} $$\Var(\Vol_k(\PVT(\gamma,W)))\sim\tau_k^{(d)}\gamma^{{d-2k\over d}}\Vol_d(W).$$ Unfortunately, no analytic expression for $\tau_k^{(d)}$ is currently known. Only for the planar case Brakke \cite{Brakke} has obtained $$\tau_1^{(2)}\approx 1.0445685$$ as the result of the evaluation of a rather involved multiple integral. That is the reason, why we will restrict in this context to the planar case $d=2$. Here we have $$J(B_R^2)={16\over 3}R^3$$ and $$\Var(\Vol_1(\PLT(t,B_R^2)))\sim{16\over 3}tR^3.$$ For the planar Poisson-Voronoi tessellation $\PVT(\gamma)$
it is well known that its edge length intensity $L_A$ is related to $\gamma$ by $L_A=2\sqrt{\gamma}$. Taking now $L_A=t$ as the parameter, we obtain $$\Var(\Vol_1(\PVT(t^2/4,B_R^2))\sim\pi\tau_1^{(2)} R^2,$$ which is in particular independent of the edge length intensity $t$. The last fact can intuitively explained as follows: per unit area we have, in mean, $t^2\over 4$ cells with mean total perimeter $c_1\cdot t$. Thus, a single cell has mean perimeter $c_2 t$ and Variance $c_3 t^2$. Since the cells are 'almost independent' of each other (see \cite[Chap. 10.5]{SW}), these variances add up and $t^2$ cancels out (here $c_1,c_2,c_3$ stand for some real constants). In fact a lot more
is known about the typical geometry of Poisson-Voronoi tessellations in the plane, including further explicit information, see e.g. Calka \cite{Ca1,Ca2}.
\begin{figure}[t]
 \begin{center}
 \includegraphics[width=10cm]{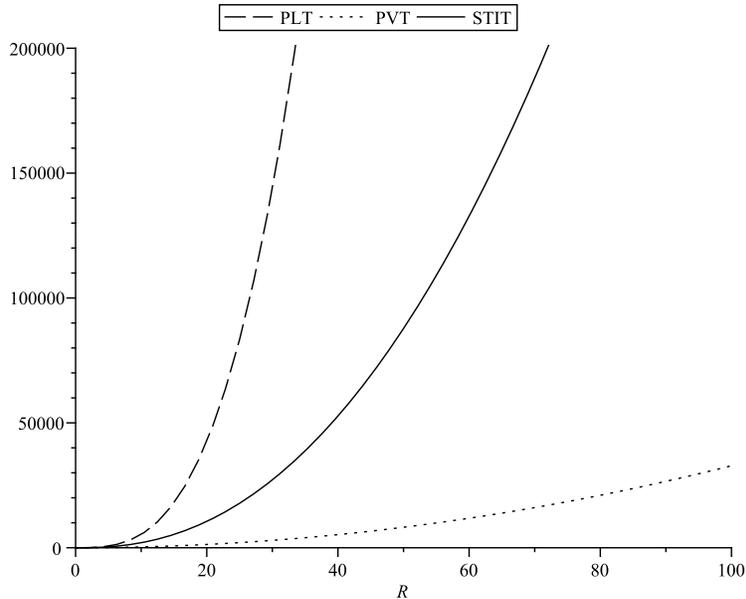}
 \caption{Comparison of the asymptotic variance of the total edge length for a STIT tessellation $Y(1)$ (thick line) with that of a Poisson line tessellation (dashed line) and a Poisson-Voronoi tessellation (thin line) with the same surface intensity in a growing ball $B_R^2$}\label{Fig3}
 \end{center}
\end{figure}
In Fig. \ref{Fig3} the graphs of the total variance of the edge length of $Y(1,B_R^2)$, $\PLT(1,B_R^2)$ and $\PVT(1/4,B_R^2)$ as a function of $R$ are compared and the behaviour discussed above is visualized. As already observed in earlier papers, STIT tessellations turn out to {\it interpolate} between Poisson line and Poisson-Voronoi tessellations in some sense. The list where this is true can now be extended by the variance of the total edge length (or more generally the total surface area).\\We can also compare the reduced second moment function and the pair-correlation function of $Y(t)$ given by (\ref{KR}) and (\ref{GR}) with the corresponding functions of a stationary and isotropic Poisson hyperplane tessellation $\PHT(t)$ with the same surface intensity $t>0$. We will denote these two functions by $K_d^{\PHT(t)}(r)$ and $g_d^{\PHT(t)}(r)$. From Slivnyak's theorem for Poisson processes one immediately infers that $tK_d^{\PHT(t)}(r)$ is given by $tK_d^{\PHT(t)}(r)=\Vol_{d-1}(B_r^{d-1})+t\Vol_d(B_r^d)$ and thus $$K_d^{\PHT(t)}(r)={\kappa_{d-1}\over t}r^{d-1}+\kappa_dr^d.$$ From the definition of the pair-correlation function we obtain now $$g_d^{\PHT(t)}(r)={1\over d\kappa_dr^{d-1}}{dK_d^{\PHT(t)}(r)\over dr}=1+{(d-1)\kappa_{d-1}\over d\kappa_dtr}.$$ 
\begin{figure}
\begin{center}
 \includegraphics[width=8cm]{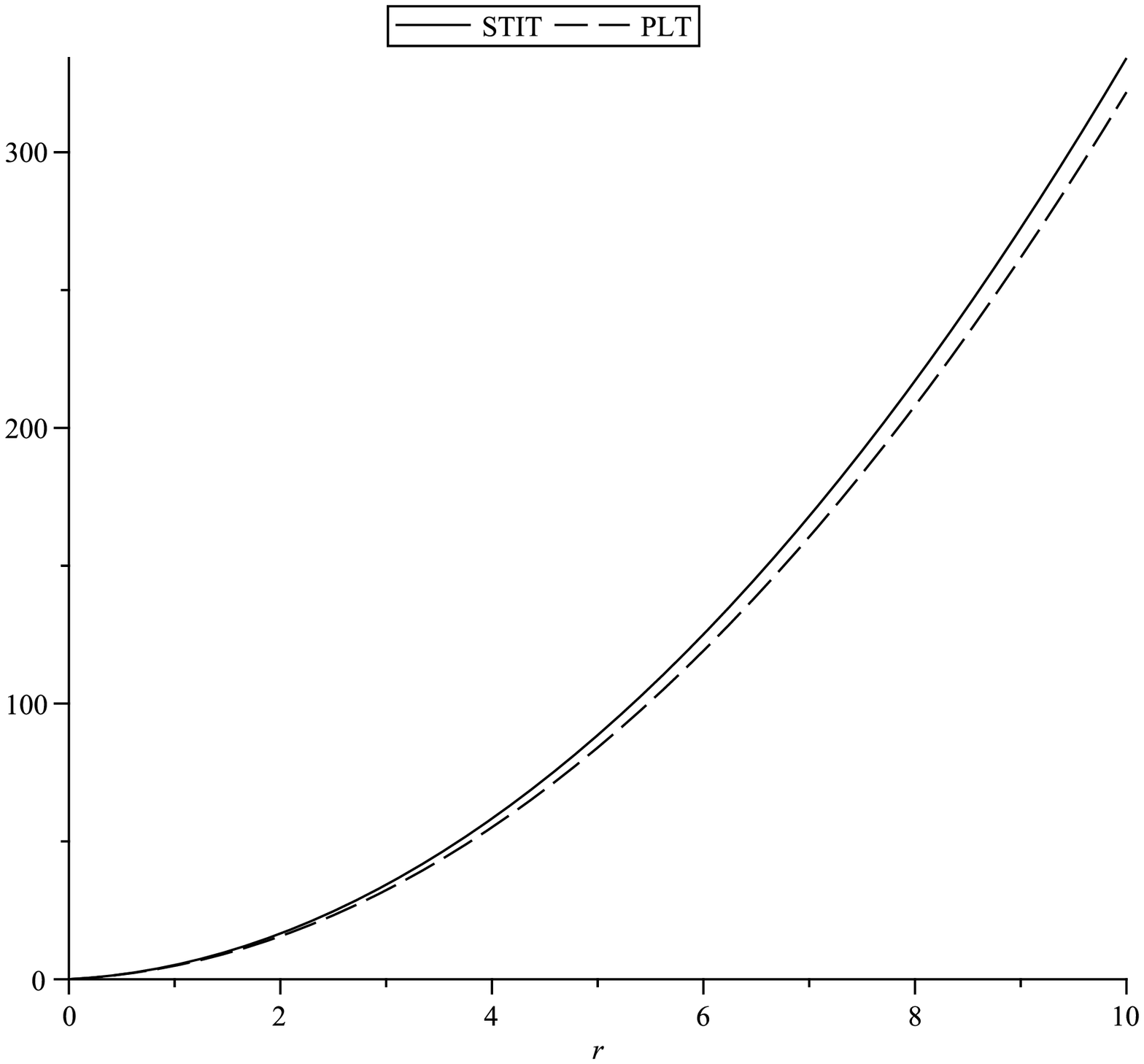}
 \includegraphics[width=8cm]{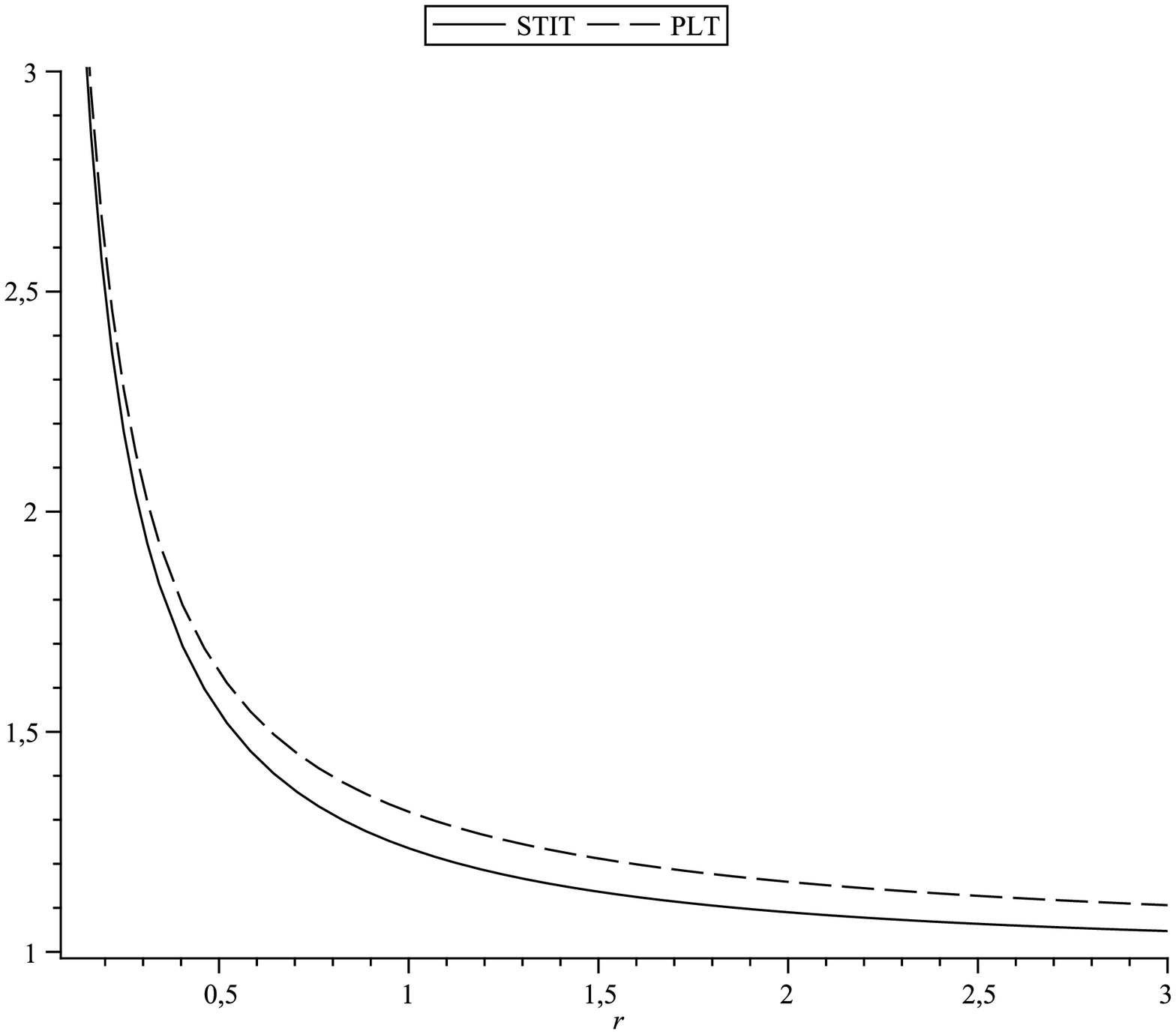}
 \caption{Fig. 3: Comparison of the $K$-function (left) and the pair-correlation function (right) for a STIT tessellation $Y(1)$ and a Poisson line tessellation $\PLT(1)$ with the same parameter}\label{Fig4}
\end{center}
\end{figure}
Especially for the planar case $d=2$, i.e. for the Poisson line tessellation $\PLT(t)$ we have the formulas $$K_2^{\PLT(t)}(r)=\pi r^2+{2r\over t}\ \ \ \text{and}\ \ \ g_2^{\PLT(t)}(r)=1+{1\over \pi tr}.$$ We have to skip the comparison of the K-function and the pair-correlation function between the STIT and Poisson-Voronoi tessellations, since nothing seems to be known about the pair-correlation function of the random length measure of a planar Poisson-Voronoi tessellation. However, the second order analysis of the point process of nodes of a Voronoi tessellations is rather involved and we refer the reader to Heinrich and Muche \cite{HM} and references therein for details, including numerics. In Fig. \ref{Fig4} we can see that the K-function and the pair-correlation function of a planar STIT tessellation and a Poisson line tessellation are very close together and it seems -- having statistical applications in mind -- that comparing variances could be more fruitful and could yield much better information about the underlying tessellation model than a comparison of the reduced second moment function or the pair correlation function (the picture is essentially the same for any space dimension).\\ In the last 15 years, central limit theorems for random tessellations have been considered by many authors, mostly for the case of stationary and isotropic Poisson hyperplane or stationary and isotropic Poisson-Voronoi tessellations, see Heinrich \cite{H} and the references therein. Elegant proofs of central limit theorems for Poisson hyperplane tessellations rely on the fact, that many functionals may be expressed as U-statistics with a special kernel function.  Application of Hoeffings's decomposition and Hoeffding's central limit theorem for U-statistics leads now to limit theorems for geometric functionals of Poisson hyperplane tessellations. The classical Cram\'er-Wold device also allows multivariate extensions. Note that in all cases the limit distribution is a normal distribution. Note that the non-additive ovoid functional $J(W)$ from (\ref{EQJD}) may also be written as $$J(W)={2\kappa_{d-1}\over d^2\kappa_d}I_d(W),$$ where $I_d$ is the $d$-th chord power integral of $W$ (here one uses again the affine Blaschke-Petkantschin formula to show the equality). It means that the asymptotic variance of the total surface area of a Poisson hyperplane tessellation depends on $W$ by $I_d(W)$ and also on the surface intensity $t$.\\ The central limit theorems for functionals of the Poisson-Voronoi tessellations, as provided by many authors \cite{BY,H,HSS,HM,Pe1,Pe2,PY1,PY2}, are mainly based on their strong mixing properties and the fact that they exhibit exponential decay of dependencies, whereas Poisson hyperplane tessellations are strongly mixing, but long-range dependent, which means that we observe slow decay of the correlations between distant parts of the tessellations. Unfortunately -- as already explained above -- the results are much less explicit in this case, since most often analytic expressions for variances are not known. But again, the limit distribution were shown to be normal distributions in all cases. The asymptotic variances are shown to depend on the shape of the sequence of growing observation windows $W_R$ only through their volume $\Vol_d(W)$ which may itself also be (artificially) expressed as ${2\over d\kappa_d}I_1(W)$, where $I_1(W)$ is the first-order chord power integral of $W$.
They also depend (beside the planar case) on the parameter $t$.\\ In the following table the facts about asymptotic variances of the total surface area of a PHT, a PVT and a STIT tessellation in the sequence of growing convex observation windows $W_R\subset{\Bbb R}^d$ are summarized for the isotropic cases ($c_d$ denotes some space dimension dependent constant, but not necessarily always the same, and $I_k(W)$ is the $k$-th chord power integral of $W$).
\begin{center}
\begin{tabular}{|c p{3.5cm}||p{2.7cm}|p{4.5cm}|p{2.7cm}|}
\hline
\parbox[0pt][2em][c]{0cm}{} & Model & PVT & STIT & PHT\\
\hline
\hline
\parbox[0pt][2em][c]{0cm}{} & $\Var(\Vol_{d-1}(\cdot,W_R))$ & $c_dt^{2-d}I_1(W)R^d$ & ${c_2I_1(W)R^2\log R\ \ (d=2)}$ ${c_dI_{d-1}(W)R^{2d-2}\ \ (d\geq 3)}$ & $c_dtI_d(W)R^{2d-1}$\\
\hline
\end{tabular}
\end{center}
There are also some central limit theorems for stationary iterated tessellations, which were motivated by potential applications in telecommunication. However, these theorems do not lead to corresponding results for STIT tessellations. As a nice survey for the central limit theory of Poisson hyperplane, Poisson-Voronoi and iterated tessellations with full proofs and several applications we recommend the small book \cite{HS}. Another 
recommendable source is the recent survey \cite{Ca3}.

\end{document}